\documentclass[a4paper, 11pt, english]{article}
\usepackage[utf8]{inputenc}
\usepackage[T1]{fontenc}
\usepackage{graphicx}
\usepackage{stmaryrd}
\usepackage[a4paper]{geometry}
\usepackage{amsmath,amsfonts,amssymb,amsthm,epsfig,epstopdf,url,array}
\usepackage{rotating}
\usepackage[colorlinks=true,citecolor=red,linkcolor=blue,pdfpagetransition=Blinds]{hyperref}
\usepackage{cleveref}
\usepackage{nameref}
\usepackage{enumitem}
\usepackage{comment}
\Crefname{paragraph}{Section}{Sections}
\usepackage{fancyhdr}

\usepackage{fullpage}

\newcommand{\ensemblenombre}[1]{\mathbb{#1}}

\newcommand{\R}{} 
\renewcommand{\R}{\ensemblenombre{R}}
\newcommand{\C}{\ensemblenombre{C}}

\newcommand{\abs}[1]{\left\lvert#1\right\rvert}
\newcommand{\norme}[1]{\left\lVert#1\right\rVert}

\newcommand{\dive}[1]{\mathrm{div}}
\newcommand{\ov}[1]{\overline{#1}}
\newcommand{\intzero}{\!\!\!\int_{\dom_0}}

\providecommand{\keywords}[1]{\noindent {\textit{Keywords:}} #1}

\theoremstyle{plain} 
\newtheorem{prop}{Proposition}[section] 
\newtheorem{theo}[prop]{Theorem}

\newtheorem{lem}[prop]{Lemma}

\theoremstyle{definition}

\newtheorem{rmk}[prop]{Remark}

\def\dx{\,\textnormal{d}x}
\def\dt{\textnormal{d}t}
\def\d{\textnormal{d}}
\def\esp{{\mathbb{E}}}
\def\dom{\mathcal{D}}

\newcommand{\vertiii}[1]{{\left\vert\kern-0.25ex\left\vert\kern-0.25ex\left\vert #1 
    \right\vert\kern-0.25ex\right\vert\kern-0.25ex\right\vert}}
\newcommand{\weakly}{\rightharpoonup}

\makeatletter
\let\original@addcontentsline\addcontentsline
\newcommand{\dummy@addcontentsline}[3]{}
\newcommand{\DeactivateToc}{\let\addcontentsline\dummy@addcontentsline}
\newcommand{\ActivateToc}{\let\addcontentsline\original@addcontentsline}
\makeatother

\begin{document}

\title{Global null-controllability for stochastic semilinear parabolic equations}
\author{V\'ictor Hern\'andez-Santamar\'ia\thanks{Instituto de Matem\'aticas, Universidad Nacional Aut\'onoma de M\'exico, Circuito Exterior, C.U., C.P. 04510 CDMX, Mexico. E-mail: \texttt{victor.santamaria@im.unam.mx}} \and K\'evin Le Balc'h \thanks{Institut de Math\'ematiques de Bordeaux, 351 Cours de la Lib\'eration, 33400 Bordeaux, France. E-mail: \texttt{kevin.le-balch@math.u-bordeaux.fr}} \and Liliana Peralta \thanks{Centro de Investigaci\'on en Matem\'aticas, UAEH, Carretera Pachuca-Tulancingo km 4.5 Pachuca, Hidalgo 42184, Mexico. E-mail: \texttt{liliana\_peralta@uaeh.edu.mx}}}

\maketitle

\begin{abstract}
In this paper, we prove the small-time global null-controllability of forward (resp. backward) semilinear stochastic parabolic equations with globally Lipschitz nonlinearities in the drift and diffusion terms (resp. in the drift term). In particular, we solve the open question posed by S. Tang and X. Zhang, in 2009. We propose a new twist on a classical strategy for controlling linear stochastic systems. By employing a new refined Carleman estimate, we obtain a controllability result in a weighted space for a linear system with source terms. The main novelty here is that the Carleman parameters are made explicit and are then used in a Banach fixed point method. This allows to circumvent the well-known problem of the lack of compactness embeddings for the solutions spaces arising in the study of controllability problems for stochastic PDEs.
\end{abstract}
\footnotesize
\tableofcontents
\normalsize
\section{Introduction}

Let $T>0$ be a positive time, $\dom$ be a bounded, connected, open subset of $\R^N$, $N\in\mathbb N^*$, with boundary $\Gamma := \partial\dom$ regular enough. Let $\dom_0$ be a nonempty open subset of $\dom$. As usual, we introduce the notation $\chi_{\dom_0}$ to refer to the characteristic function of the set $\dom_0$. To abridge the notation, hereinafter we write $Q_T:=(0,T)\times \dom$ and $\Sigma_T:=(0,T)\times \Gamma$. 

Let $(\Omega, \mathcal{F}, \{\mathcal{F}_t\}_{t \geq 0}, \mathbb{P})$ be a complete filtered probability space on which a one-dimensional standard Brownian motion $\{ W(t)\}_{t \geq 0}$ is defined such that $\{\mathcal{F}_t\}_{t \geq 0}$ is the natural filtration generated by $W(\cdot)$ augmented by all the $\mathbb{P}$-null sets in $\mathcal{F}$. Let $X$ be a real Banach space, for every $p \in [1, +\infty]$, we introduce
\begin{equation*}
L_{\mathcal{F}}^p(0,T;X) := \{ \phi : \phi\ \text{is an }X\text{-valued } \mathcal{F}_t\text{-adapted process on } [0,T]\ \text{and}\ \phi \in L^p([0,T] \times \Omega ;X)\},
\end{equation*}
endowed with the canonical norm and we denote by $L_{\mathcal{F}}^2(\Omega; C([0,T];X))$ the Banach space consisting on all $X$-valued $\mathcal{F}_t$-adapted process $\phi(\cdot)$ such that $\mathbb{E}\left(\norme{\phi(\cdot)}_{C([0,T];X)}^2 \right) < \infty$, also equipped with the canonical norm.

Let us consider the stochastic forward semilinear equation 
\begin{equation}\label{eq:forward_semilinear}
\begin{cases}
\d{y}=(\Delta y + f(\omega,t,x,y)+\chi_{\dom_0}h)\dt+(g(\omega,t,x,y)+H)\d W(t) &\text{in }Q_T, \\
y=0 &\text{on }\Sigma_T, \\
y(0)=y_0 &\text{in }\dom.
\end{cases}
\end{equation}

In \eqref{eq:forward_semilinear}, $y$ is the state variable, the pair $(h,H)\in L^2_{\mathcal F}(0,T;L^2(\dom_0))\times L^2_{\mathcal F}(0,T;L^2(\dom))$ is the control and $y_0\in L^2(\Omega,\mathcal F_0; L^2(\dom))$ is the initial datum. Here, $f$ and $g$ are two globally Lipschitz nonlinear functions, that is, there exists a positive constant $L>0$ such that
\begin{align}
\begin{cases}
&\forall (\omega,t,x, s_1,s_2)\in \Omega\times[0,T]\times \dom\times \R^2,\label{eq:UniformLipschitzf_forw}\\
& |f(\omega,t,x,s_1)-f(\omega,t,x,s_2)| + |g(\omega,t,x,s_1)-g(\omega,t,x,s_2)| \leq L|s_1-s_2|.
\end{cases}
\end{align}
We also impose 
\begin{equation}
\label{eq:f_g_zero}
\forall (\omega,t, x) \in \Omega\times[0,T]\times\dom,\  f(\omega,t,x,0) = g(\omega,t,x,0)=0.
\end{equation}

Under these conditions, by taking $y_0\in L^2(\Omega,\mathcal F_0; L^2(\dom))$ and $(h,H)\in L^2_{\mathcal F}(0,T;L^2(\dom_0))\times L^2_{\mathcal F}(0,T;L^2(\dom))$, it is known (see \cite[Theorem 2.7]{LZ19}) that system \eqref{eq:forward_semilinear} is globally defined in $[0,T]$. More precisely, we can establish the existence and uniqueness of the solutions to \eqref{eq:forward_semilinear} in the class 
\begin{equation}\label{eq:wp_forward}
y\in \mathcal W_T:= L^2_{\mathcal F}(\Omega; C([0,T];L^2(\dom)))\cap L^2_{\mathcal F}(0,T; H_0^1(\dom)).
\end{equation}

One of the key questions in control theory is to determine whether a system enjoys the so-called null controllability property. System \eqref{eq:forward_semilinear} is said to be globally null-controllable if for any initial datum $y_0\in L^2(\Omega,\mathcal F_0; L^2(\dom))$, there exist controls $(h,H)\in L^2_{\mathcal F}(0,T;\dom_0)\times L^2_{\mathcal F}(0,T;L^2(\dom))$ such that the corresponding solution satisfies
\begin{equation}\label{eq:null_cond}
y(T,\cdot)=0 \quad\text{in $\dom$, a.s.}
\end{equation}
Observe that the regularity \eqref{eq:wp_forward} justifies the definition we have introduced. 

In this paper, we are interested in studying this controllability notion for system \eqref{eq:forward_semilinear}.  Before introducing our main results we give a brief panorama of previous results available in the literature and emphasize the main novelty of this work. 

\subsection{Known results}

The controllability of parabolic partial differential equations (PDEs) has been studied by many authors and the results available in the literature are very rich. In the following paragraphs, we focus on (small-time) global null-controllability results for scalar parabolic equations.

\textbf{Deterministic setting.} In the case where $g\equiv H\equiv 0$ and $f$ and $y_0$ are deterministic functions, system \eqref{eq:forward_semilinear} has been studied by several authors. In the mid 90's, Fabre, Puel \& Zuazua in \cite{FPZ95} studied the so-called global approximate-null controllability in the case where $f$ is a globally Lipschitz nonlinearity  and condition \eqref{eq:null_cond} is replaced by the weaker constraint $\|y(T)\|_{L^2(\dom)}\leq \epsilon$. Later, Imanuvilov \cite{Ima95} and Fursikov \& Imanuvilov \cite{fursi} improved this result and proved that the global null-controllability holds, see also \cite{LR95} for the case of the (linear) heat equation, i.e. $f \equiv 0$. After these seminal works, Fern\'andez-Cara \cite{FC97}, Fern\'andez-Cara \& Zuazua \cite{FCZ00} have considered slightly superlinear functions $f$ leading to blow-up without control, see also \cite{Bar00} and the more recent work \cite{LB20}. Results for nonlinearities including $\nabla y$ and depending on Robin boundary conditions have also been studied for instance in \cite{DFCGBZ02,FCGBGP06}.

One common feature among these results is that the authors study the controllability problem by using the following general strategy, due to Zuazua in the context of the wave equation, see \cite{Zua93} or \cite[Chapter 4.3]{Cor07} for a presentation of this method. First, linearize the system and study the controllability of the system \eqref{eq:forward_semilinear} replacing $f$ by $a(t,x) y$ where $a \in L^{\infty}(Q_T)$. Then, use a suitable fixed point method (commonly Schauder or Kakutani) for addressing the controllability of the nonlinear system. At this point, the important property of compactness is needed. In fact, compact embeddings relying on Aubin-Lions lemma like $W(0,T):=\{y\in L^2(0,T;H_0^1(\dom)),\ y_t\in L^2(0,T;H^{-1}(\dom))\} \hookrightarrow L^2(0,T;L^2(\dom))$ are systematically used. 

\textbf{Stochastic setting.} In the case where $f(y)=\alpha y$ and $g(y)=\beta y$, $\alpha,\beta\in \R$, the controllability results for \eqref{eq:forward_semilinear} were initiated by Barbu, R\u{a}\c{s}canu \& Tessitore in \cite{BRT03}. Under some restrictive conditions and without introducing the control $H$ on the diffusion, they established a controllability result for linear forward stochastic PDEs. Later, Tang \& Zhang in \cite{TZ09} improved this result and considered more general coefficients $\alpha$ and $\beta$ (depending on $t,x$ and $\omega$). The main novelty in that work was to introduce the additional control $H$ and prove fine Carleman estimates for stochastic parabolic operators. The same methodology has been used to study other cases like the ones of Neumann and Fourier boundary conditions (\cite{Yan18}), degenerate equations (\cite{LY19}) and fourth-order parabolic equations (\cite{GCL15}). As a side note, we shall mention the work by L\"{u} in \cite{LU11} who, by using the classical Lebeau-Robbiano strategy (\cite{LR95}), noticed that the action of the control $H$ can be omitted at the prize of considering random coefficients $\alpha$ and $\beta$ only depending on the time variable $t$. 

In the framework proposed in this paper, as far the author's knowledge, there are not any results available in the literature. Compared to the deterministic setting, while establishing controllability properties for stochastic PDEs, many new difficulties arise. For instance, the solution of stochastic PDEs are usually not differentiable with respect to the variable with noise (i.e., the time variable). Also the diffusion term introduces additional difficulties while analyzing the problem. But most importantly, as remarked in \cite[Remark 2.5]{TZ09}, the compactness property, which is one of the key tools in the deterministic setting, is known to be false for the functional spaces related to stochastic PDEs. This is the main obstruction for employing some classical methodologies like in \cite{fursi}, \cite{FCZ00} for establishing null-controllability of semilinear problems at the stochastic level.

\subsection{Statement of the main results}

To overcome the lack of compactness mentioned in the last section, in this work we propose a new tweak on an old strategy for controlling parabolic systems. We use a classical methodology for controlling a linear system with source terms $F,G\in L^2_{\mathcal F}(0,T;L^2(\dom))$ of the form
\begin{equation}\label{eq:sys_forward_source_intro}
\begin{cases}
\d{y}=(\Delta y + F+ \chi_{\dom_0}h)\dt + (G+H)\d{W}(t) &\text{in }Q_T, \\
y=0 &\text{on } \Sigma_T, \\
y(0)=y_0 &\text{in }\dom,
\end{cases}
\end{equation}
in a suitable weighted space. Note that this strategy has been widely used in the literature and it has been revisited in \cite{LLT13} to obtain local results. In turn, such weighted space is naturally defined through the weights arising in the Carleman estimates needed for studying the observability of the corresponding linear adjoint system. But, unlike many other works out there, we make precise the dependency on the parameters involved in the construction of the Carleman weights and use them in a second stage to prove that the nonlinear map $\mathcal N(F,G)\mapsto (f(t,x,\omega,y),g(t,x,\omega,y))$, with $y$ solution of \eqref{eq:sys_forward_source_intro}, is well-defined and is strictly contractive in a suitable functional space. In this way, the controllability of \eqref{eq:forward_semilinear} is ensured through a Banach fixed point method which does not rely on any compactness argument. 

Our first main result is as follows. 
\begin{theo}
\label{th:semilinear_forward}
Under assumptions \eqref{eq:UniformLipschitzf_forw}--\eqref{eq:f_g_zero}, system \eqref{eq:forward_semilinear} is small-time globally null-controllable, i.e. for every $T>0$ and for every $y_0 \in L^2(\Omega,\mathcal F_0;L^2(\dom))$, there exists $h \in L^2_{\mathcal F}(0,T;L^2(\dom_0))$ such that the unique solution $y$ of \eqref{eq:forward_semilinear} satisfies $y(T,\cdot) = 0$ in $\dom$, a.s.
\end{theo}

\begin{rmk}
As compared to some results in the deterministic framework, on the one hand notice that here we are not considering any differentiability condition on the nonlinearities, that is, $f,g$ are merely  $C^0$-functions. On the other hand, our method does not permit to establish global controllability result for slightly superlinearities as considered in \cite{FCZ00}.
\end{rmk}

\begin{rmk}\label{rmk:carleman_new}
As mentioned above, the main ingredient to prove \Cref{th:semilinear_forward} is a precise Carleman estimate for the linear adjoint system to \eqref{eq:sys_forward_source_intro} (see \Cref{thm:carleman_backward_0}), which in this case is a backward parabolic equation. Previous to this work, such Carleman estimate was not available in the literature (see \cite[Theorem 2.5]{BEG16} for a similar estimate in the deterministic case). The methodology employed to prove the result is the weighted identity method introduced in the stochastic framework in \cite{TZ09}. 
\end{rmk}

As classical in the stochastic setting, for completeness, using the same strategy described above, it is possible to establish a controllability result for semilinear backward parabolic equations. More precisely, consider 

\begin{equation}\label{eq:backward_semilinear}
\begin{cases}
\d{y}=\left(-\Delta y+f(\omega,t,x,y,Y)+\chi_{\dom_0}h\right)\d{t}+Y\d{W}(t) &\text{in } Q_T, \\
y=0 &\text{on }\Sigma_T, \\
y(T)=y_T &\text{in }\dom,
\end{cases}
\end{equation}
where $f$ is a globally Lipschitz nonlinearity, i.e., there exists $L>0$ such that
\begin{align}
\label{eq:UniformLipschitzf}
\begin{cases}
\forall (w,t, x) \in \Omega\times [0,T]\times\dom,\ \forall s_1,s_2,\ov{s}_1,\ov{s}_2 \in \R, \\ |f(\omega,t,x,s_1,\ov{s}_1) - f(\omega,t,x,s_2,\ov{s}_2)| \leq L \left( |s_1-s_2|+ |\ov{s}_1-\ov{s}_2|\right).
\end{cases}
\end{align}
Moreover, we impose that 
\begin{equation}
\label{eq:fzero}
\forall (\omega,t, x) \in \Omega\times[0,T]\times\dom,\  f(\omega,t,x,0,0) = 0.
\end{equation}

Under these conditions, by taking $y_T\in L^2(\Omega,\mathcal F_T; L^2(\dom))$ and $h \in L^2_{\mathcal F}(0,T;L^2(\dom_0))$, it is known (see \cite[Theorem 2.12]{LZ19}) that system \eqref{eq:backward_semilinear} is also globally well-defined in $[0,T]$. In this case, we can establish the existence and uniqueness of the solutions to \eqref{eq:backward_semilinear} in the class 
\begin{equation}\label{eq:wp_backward}
(y,Y)\in \mathcal W_T\times L^2_{\mathcal F}(0,T;L^2(\dom)).
\end{equation}

Our second main result is as follows.
\begin{theo}
\label{th:semilinear_backward}
Under assumptions \eqref{eq:UniformLipschitzf}--\eqref{eq:fzero}, system \eqref{eq:backward_semilinear} is small-time globally null-controllable, i.e. for every $T>0$ and for every $y_T \in L^2(\Omega,\mathcal F_T;L^2(\dom))$, there exists $h \in L^2_{\mathcal F}(0,T;L^2(\dom_0))$ such that the unique solution $y$ of \eqref{eq:backward_semilinear} satisfies $y(0,\cdot) = 0$ in $\dom$, a.s.
\end{theo}

\Cref{th:semilinear_backward} extends to the nonlinear setting the previous results in
\cite[Corollary 3.4]{BRT03} and \cite[Theorem 2.2]{TZ09} for the backward equation.

The strategy to prove \Cref{th:semilinear_backward} is very close to the one of \Cref{th:semilinear_forward} but one major difference can be spotted. For this case, it is not necessary to prove a Carleman estimate for forward stochastic parabolic equations. Actually, it suffices to use the deterministic Carleman inequality of \cite[Thm. 2.5]{BEG16} and employ the duality method introduced by Liu in \cite{liu14}.

\subsection{Outline of the paper}

The rest of the paper is organized as follows. In \Cref{sec:forward}, we present the proof of \Cref{th:semilinear_forward}. In particular, \Cref{sec:new_carleman} is devoted to prove the new Carleman estimate anticipated in \Cref{rmk:carleman_new}. \Cref{sec:backward} is devoted to prove \Cref{th:semilinear_backward}. Finally in \Cref{sec:conclusion} we present some concluding remarks.

\section{Controllability of a semilinear forward stochastic parabolic equation}\label{sec:forward}

\subsection{A new global Carleman estimate for a backward stochastic parabolic equation}\label{sec:new_carleman}

This section is devoted to prove a new Carleman estimate for a backward stochastic parabolic equation. The main novelty here is that the weight does not degenerate as $t\to 0^+$ (compared with the classical work \cite{fursi}). This estimate has been proved in the deterministic case in \cite[Theorem 2.5]{BEG16} in a slightly more general framework. Here, we use many of the ideas presented there and adapt them to the stochastic setting.

To make a precise statement of our result, let $\dom^\prime$ be a nonempty subset of $\dom$ such that $\dom^\prime\subset\subset\dom_0$. Let us introduce $\beta\in C^4(\ov{\dom})$ such that
\begin{equation}\label{eq:prop_weight_beta}
\begin{cases}
0<\beta(x)\leq 1, \ \forall x\in{\dom}, \ \\
\beta(x)=0, \ \forall x\in \partial\dom, \\
\inf_{\dom\setminus\ov{\dom^\prime}} \{|\nabla \beta|\}\geq \alpha >0.
\end{cases}
\end{equation}
The existence of such a function is guaranteed by \cite[Lemma 1.1]{fursi}. 

Without loss of generality, in what follows we assume that $0<T<1$. For some constants $m\geq 1$ and $\sigma\geq 2$ we define the following weight function depending on the time variable

\begin{equation}\label{eq:def_theta}
\gamma(t):=
\begin{cases}
\gamma(t)=1+(1-\frac{4t}{T})^\sigma, \quad t\in(0,T/4], \\
\gamma(t)=1, \quad t\in[T/4,T/2], \\
\gamma \textnormal{ is increasing on $[T/4,T/2]$}, \\
\gamma(t)=\frac{1}{(T-t)^m}, \quad t\in[3T/4,T], \\
\gamma\in C^2([0,T)).
\end{cases}
\end{equation}
We take the following weight functions $\varphi=\varphi(t,x)$ and $\xi=\xi(t,x)$
\begin{equation}\label{eq:weights_0}
\varphi(t,x):=\gamma(t)\left(e^{\mu(\beta(x)+6m)}-\mu e^{6\mu(m+1)}\right), \quad \xi(t,x):=\gamma(t)e^{\mu(\beta(x)+6m)},
\end{equation}
where $\mu$ is a positive parameter with $\mu\geq 1$ and $\sigma$ is chosen as
\begin{equation}\label{def:sigma}
\sigma=\lambda\mu^2e^{\mu(6m-4)},
\end{equation}
for some parameter $\lambda\geq 1$. Observe that with these elections on $\mu$ and $\lambda$, the parameter $\sigma$ is always greater than 2 and this also ensures that $\gamma(t)\in C^2([0,T))$. We finally set the weight $\theta=\theta(t,x)$ as 
\begin{equation}\label{def:theta_ell}
\theta:=e^{\ell} \quad\text{where } \ell(t,x):=\lambda\varphi(t,x).
\end{equation}

Using these notations, we state the main result of this section which is a Carleman estimate for backward stochastic parabolic equations.

\begin{theo} \label{thm:carleman_backward_0}
For all $m\geq 1$, there exist constants $C>0$, $\lambda_0\geq 1$ and $\mu_0\geq 1$ such that, for any $z_T\in L^2(\Omega,\mathcal F_T;\dom)$ and any $\Xi\in L^2_{\mathcal F}(0,T;L^2(\dom))$, the solution $(z,\ov{z})\in \mathcal W_T \times L^2_{\mathcal F}(0,T;L^2(\dom))$ to
\begin{equation}\label{eq:system_z}
\begin{cases}
\d{z}=(-\Delta z+\Xi)\dt+\ov{z}\d{W}(t) &\text{in }Q_T, \\
z=0 &\text{on }\Sigma_T, \\
z(T)=z_T &\text{in }\dom,
\end{cases}
\end{equation}
satisfies
\begin{align}\notag
\esp&\left(\int_{\dom}\theta^2(0)|\nabla z(0)|^2\dx\right)+\esp\left(\int_{\dom}\lambda^2\mu^3e^{2\mu(6m+1)}\theta^2(0)|z(0)|^2\dx\right)\\ \notag
&+\esp\left(\int_{Q_T}\lambda\mu^2\xi\theta^2|\nabla z|^2\dx\dt\right)+\esp\left(\int_{Q_T}\lambda^3\mu^4\xi^3\theta^2| z|^2\dx\dt\right) \\ \label{eq:car_backward_0}
&\leq C\esp\left(\int_{0}^{T}\!\!\!\int_{\dom_0}\lambda^3\mu^4\xi^3\theta^2|z|^2\dx\dt+\int_{Q_T}\theta^2|\Xi|^2\dx\dt+\int_{Q_T}\lambda^2\mu^2\xi^3\theta^2|\ov{z}|^2\dx\dt\right)
\end{align}
for all $\mu\geq \mu_0$ and $\lambda\geq \lambda_0$. 
\end{theo}

Before giving the proof of \Cref{thm:carleman_backward_0}, we make the following remark.

\begin{rmk}
We make the following comments.
\begin{itemize}
\item The proof of this result is rather classical except for the definition of the weight function $\gamma(t)$ which does not blow-up as $t\to 0^+$ and thus preventing $\theta$ to vanish at $t=0$. This change introduces some additional difficulties to the classical proof of the Carleman estimate for backward stochastic parabolic equations shown in \cite[Theorem 6.1]{TZ09}, but which can be handled just as in the deterministic case (see \cite[Appendix A.1]{BEG16}). 
\item Different to the Carleman estimate in \cite[Eq. (6.2)]{TZ09}, the power of $\xi$ in the last term of  \eqref{eq:car_backward_0} is $3$ rather than $2$. This is due to the definition of the weight $\gamma$ which modifies a little bit the estimate of $\varphi_t$ in $[0,T/4]$ as compared, for instance, to \cite{TZ09}. This does not represent any problem for proving our main controllability result. 
\item Just as in \cite[Remark 6.1]{TZ09}, we can estimate the last term in \eqref{eq:car_backward_0} by  weighted integrals of $f$ and $z$, more precisely
\begin{equation*}
\esp\left(\int_{Q_T}\lambda^2\mu^2\xi^3\theta^2|\ov{z}|^2\dx\dt\right)\leq C\esp\left(\int_{Q_T}\lambda^4\mu^4\xi^6\theta^2|z|^2\dx\dt+\int_{Q_T}\theta^2|\Xi|^2\dx\dt\right).
\end{equation*}
Nevertheless, the new term of $z$ cannot be controlled by its counterpart in the left-hand side of \eqref{eq:car_backward_0} and this does not improve our result. 
\end{itemize}
\end{rmk}

\begin{proof}[Proof of \Cref{thm:carleman_backward_0}]
As we have mentioned before, the proof of this result is close to other proofs for Carleman estimates in the stochastic setting (see, e.g. \cite{TZ09,Yan18} or \cite[Ch. 3]{FLZ19}). Some of the estimates presented in such works are valid in our case but others need to be adapted. For readability, we have divided the proof in several steps and we will emphasize the main changes with respect to previous works. 

\smallskip
\textbf{Step 1. A point-wise identity for a stochastic parabolic operator}. We set $\theta=e^{\ell}$ where we recall that $\ell=\lambda\varphi$ with $\varphi$ is defined in \eqref{eq:weights_0}. Then, we write $\psi=\theta z$ and for the operator $\d z+\Delta z \dt$ we have the following identity
\begin{equation}\label{eq:iden_operator}
\theta(\d{z}+\Delta z \dt)=I_1+I\dt,
\end{equation}
where
\begin{align}\label{defs:init}
\begin{cases}
I_1=\d\psi-2\sum_{i}\ell_i\psi_i\dt+\Psi\psi\dt, \\
I= A\psi+\sum_{i}\psi_{ii}, \\
A=-\ell_t+\sum_{i}(\ell_i^2-\ell_{ii})-\Psi,
\end{cases}
\end{align}
where $\Psi=\Psi(x,t)$ is a function to be chosen later. Hereinafter, to abridge the notation, we simply write $\rho_i=\partial_{x_i}\rho$ and $\rho_t=\partial_t\rho$ and we denote $\sum_i$ and $\sum_{i,j}$ to refer to $\sum_{i=1}^{N}$ and $\sum_{i=1}^{N}\sum_{j=1}^{N}$, respectively.

From It\^{o}'s formula, we have that
\begin{align*}
\d(A\psi^2)&=A\psi\d{\psi}+\psi\d(A\psi)+\d(A\psi)\d{\psi} \\ 
&= 2A\psi\d{\psi}+A_t\psi^2\d{\psi}+2\psi\d{A}\d{\psi}+A(\d{\psi})^2+(\d{A})(\d{\psi})^2
\end{align*}
and 
\begin{align*}
\psi_{ii}\d{\psi}=(\psi_i\d{\psi})_i-\psi_i\d{\psi_i}=(\psi_i\d{\psi})_i-\frac{1}{2}\d(\psi_i^2)+\frac{1}{2}(\d\psi_i)^2.
\end{align*}
Therefore
\begin{align}\notag 
I\d{\psi}&=\left(A\psi+\sum_{i}\psi_{ii}\right) \d \psi \\ \label{eq:iden_I_dpsi}
&= \sum_{i} (\psi_i \d\psi)_i-\frac{1}{2}\d\left(\sum_{i}\psi_i^2\right)+\frac{1}{2}\sum_{i}(\d\psi_i)^2+\frac{1}{2}\d(A\psi^2)-\frac{1}{2}A_t\psi^2\dt-\frac{1}{2}A(\d{\psi})^2.
\end{align}
On the other hand, a direct computation gives
\begin{align} \notag
-2\sum_i\ell_i\psi_i I=&-\sum_i(A\ell_i\psi^2)_i+\sum_i(A\ell_i)_i\psi^2 \\ \label{eq:iden_li_psii_I}
&-\sum_i\Big[\sum_j(2\ell_j\psi_i\psi_j-\ell_i\psi_j\psi_j)\Big]_i+\sum_{i,j}\sum_{k,h}\left[2\delta_{ih}\delta_{kj}\ell_{kh}-\delta_{ij}\delta_{kh}\ell_{kh}\right]\psi_i\psi_j,
\end{align}
where $\delta_{ij}=1$ if $i=j$ and 0 otherwise. Multiplying both sides of \eqref{eq:iden_operator} by $I$ and taking into account identities \eqref{eq:iden_I_dpsi}--\eqref{eq:iden_li_psii_I}, we get the following point-wise identity 
\begin{align}\notag
\theta I(\d z+\Delta z\dt)&=\left(I^2+\sum_{i,j}c^{ij}\psi_i\psi_j+F \psi^2+I\Psi\psi+\nabla\cdot V \right)\dt+\sum_{i}(\psi_i \d\psi)_i \\ \notag
&\quad +\frac{1}{2}\sum_{i}(\d\psi_i)^2-\frac{1}{2}A(\d\psi)^2 \\ \label{eq:pointwise_iden}
&\quad -\frac{1}{2}\d\left(\sum_{i}\psi_i^2-A\psi^2\right),
\end{align}
where
\begin{align*}
\begin{cases}
V=\left[V^1,V^2,\ldots,V^N\right], \\
V^i=-\sum_{j}\left(2\ell_j\psi_i\psi_j-\ell_i\psi_j\psi_j\right)-A\ell_i\psi^2, \quad i=1,\ldots,N, \\
c^{ij}=\sum_{k,h}\left[2\delta_{ih}\delta_{kj}\ell_{kh}-\delta_{ij}\delta_{kh}\ell_{kh}\right], \\
F=-\frac{1}{2}A_t\psi^2+\sum_{i}(A\ell_i)_i\psi^2.
\end{cases}
\end{align*}

\smallskip
\textbf{Step 2. Some old and new estimates}. The main goal of this step is to start building our Carleman estimate taking as a basis the point-wise identity \eqref{eq:pointwise_iden}. Integrating with respect to time in both sides of \eqref{eq:pointwise_iden}, we get
\begin{align}\label{eq:point_iden_lhs}
\int_0^{T}&\theta I(\d z+\Delta z\dt) \\ \label{eq:0_T_term}
&= -\frac{1}{2}\left.\left(\sum_{i}\psi_i^2-A\psi^2\right)\right|_{0}^{T} \\ \label{eq:F_term}
&\quad + \int_{0}^{T}\left(F \psi^2+I\Psi\psi\right)\dt \\ \label{eq:comb_terms}
&\quad +\int_{0}^{T}I^2\dt+\int_0^{T}\sum_{i,j}c^{ij}\psi_i\psi_j\dt+\int_0^{T}\nabla\cdot V \dt+\int_{0}^{T}\sum_{i}(\psi_i \d\psi)_i \\ \label{eq:point_last_rhs}
&\quad +\int_{0}^{T}\frac{1}{2}\sum_{i}(\d\psi_i)^2-\int_{0}^{T}\frac{1}{2}A(\d\psi)^2,
\end{align}
for a.e. $x\in \R^N$ and a.s. $\omega\in \Omega$. We will pay special attention to terms \eqref{eq:0_T_term} and \eqref{eq:F_term} which yield positive terms that are not present in other Carleman estimates using the classical weight vanishing both at $t=0$ and $t=T$. 

At this point, we shall choose the function $\Psi=\Psi(x,t)$ as 
\begin{equation}\label{eq:definicion_Psi}
\Psi:=-2\sum_{i}\ell_{ii}
\end{equation}
and, for convenience, we give some identities that will be useful in the remainder of the proof. From the definition of $\ell$ in \eqref{def:theta_ell}, we have
\begin{equation}\label{eq:deriv_weights}
\begin{split}
\ell_i&=\lambda\gamma\mu\beta_ie^{\mu(\beta+6m)}, \\
\ell_{ii}&=\lambda\gamma\mu^2\beta_i^2e^{\mu(\beta+6m)}+\lambda\gamma\mu\beta_{ii}e^{\mu(\beta+6m)}.
\end{split}
\end{equation}
For shortness, we have dropped the explicit dependence of $x$ and $t$ on the expressions above.

\textit{-Positivity of the term \eqref{eq:0_T_term}}. From the definitions of $\psi$ and $\ell$, we readily see that $\lim_{t\to T^{-}}\ell(t,\cdot)=-\infty$ and thus the term at $t=T$ vanishes. Therefore, \eqref{eq:0_T_term} simplifies to 
\begin{equation}\label{eq:term_at_0}
-\frac{1}{2}\left.\left(\sum_{i}\psi_i^2-A\psi^2\right)\right|_{0}^{T}=\frac{1}{2}\left(\sum_{i}\psi_i^2(0)+A(0)\psi^2(0)\right).
\end{equation}
It is clear that the first term in the right-hand side of \eqref{eq:term_at_0} is positive. For the second one, we will generate a positive term by using the explicit expression of the function $\gamma$. Using definition \eqref{eq:def_theta}, we obtain that
\begin{equation*}
\gamma^\prime(t)=-\frac{4\sigma}{T}\left(1-\frac{4t}{T}\right)^{\sigma-1}, \quad \forall t\in[0,T/4],
\end{equation*}
whence, from \eqref{def:sigma} and the above expression, we get
\begin{align}\notag 
\ell_t(0,\cdot)&=-\frac{4\lambda^2\mu^2e^{\mu(6m-4)}}{T}\left(e^{\mu(\beta(\cdot)+6m)}-\mu e^{6\mu(m+1)}\right) 
\\ \label{eq:est_lt_0}
 &\geq  c\lambda^2\mu^3 e^{\mu(12m+2)}
\end{align}
for all $\mu\geq 1$ and some constant $c>0$ uniform with respect to $T$. On the other hand, from the derivatives \eqref{eq:deriv_weights} and using the facts that $\gamma(0)=2$ and $\beta\in C^4(\ov{\dom})$, we get
\begin{equation}\label{eq:est_A_0}
|\ell_i^2(0,\cdot)+\ell_{ii}(0,\cdot)|\leq C\lambda^2\mu^2e^{2\mu(6m+1)}.
\end{equation}

In this way, by using \eqref{eq:term_at_0}, the definition of $A$ in \eqref{defs:init} and estimates \eqref{eq:est_lt_0}--\eqref{eq:est_A_0}, there exists $\mu_1>0$, such that for all $\mu\geq \mu_1\geq 1$, we get
\begin{align}\notag
-\frac{1}{2}\left.\left(\sum_{i}\psi_i^2-A\psi^2\right)\right|_{0}^{T}&\geq \frac{1}{2}|\nabla \psi(0)|^2+c\lambda^2\mu^3e^{2\mu(6m+1)}\psi^2(0)-C\lambda^2\mu^2e^{2\mu(6m+1)}\psi^2(0) \\ \label{eq:est_final_t0}
&\geq c_1|\nabla\psi(0)|^2+c_1\lambda^2\mu^3e^{2\mu(6m+1)}\psi^2(0),
\end{align}
for some constant $c_1> 0$ only depending on $\dom$ and $\dom^\prime$.

\textit{-Estimate of the term \eqref{eq:F_term}}. This term is the most cumbersome one since the combination of some terms of $F\psi^2$ and $I\Psi\psi$ will yield a positive term that does not appear in the classical Carleman estimate with weight vanishing at $t=0$ and $t=T$.

Recalling the definition of $A$, we see that that the first term in \eqref{eq:F_term} can be written as
\begin{equation}\label{eq:F_extended}
\int_{0}^{T}F\psi^2\dt=\int_{0}^{T}\left(F_1+F_2+F_3\right)\psi^2\dt,
\end{equation}
where
\begin{align}\label{eq:defs_F}
F_1=\frac{1}{2}\ell_{tt}, \quad 
F_2=-\frac{1}{2}\sum_{i}(\ell_i^2+\ell_{ii})_t, \quad 
F_3=\sum_{i}(A_i\ell_i+A\ell_{ii}).
\end{align}

For the first term of \eqref{eq:F_extended} we argue as follows. For $t\in(0,T/4)$, using the definition of $\gamma(t)$, it is not difficult to see that $|\gamma_{tt}|\leq C\lambda^2\mu^4e^{2\mu(6m-4)}$
thus 
\begin{equation}\label{eq:est_ltt_0_T4}
|\ell_{tt}|\leq C\lambda^3\mu^5e^{2\mu(6m-4)}e^{6\mu(m+1)}\leq C\lambda^3\mu^2\xi^3
\end{equation}
where we recall that $\xi=\xi(t,x)$ is defined in \eqref{eq:weights_0}. Here, we also have used that $\mu^3e^{-2\mu}<1/2$ for all $\mu>1$. 

For $t\in(T/2,T)$, using once again the definition of $\gamma$ we have $|\gamma_{tt}|\leq C\gamma^3$.
Noting that $\varphi_{tt}=\frac{\gamma_{tt}}{\gamma}\varphi$ and using the estimate $|\varphi\gamma|\leq \mu\xi^2$ we get 
\begin{equation}\label{eq:est_ltt_T2_T}
|\ell_{tt}|=\left|\lambda\frac{\gamma_{tt}}{\gamma}\varphi\right| \leq C\lambda \mu\xi^3.
\end{equation}
Since obviously $\ell_{tt}$ vanishes for $t\in(T/4,T/2)$, we can put together estimates \eqref{eq:est_ltt_0_T4} and \eqref{eq:est_ltt_T2_T} to deduce that
\begin{equation}\label{eq:est_F1}
\int_{0}^{T}F_1\psi^2\dt\geq -C\lambda^3\mu^2\int_{0}^{T}\xi^3\psi^2\dt.
\end{equation}

We move now to the second and third terms of \eqref{eq:F_extended}. To abridge the notation, in what follows, we set
\begin{equation*}
\alpha(x):=e^{\mu(\beta(x)+6m)}-\mu e^{6\mu(m+1)}.
\end{equation*}
Notice that $\alpha(x)<0$ for all $x\in\dom$. 

From \eqref{eq:deriv_weights}, a direct computation yields
\begin{align}\notag
(\ell_i^2+\ell_{ii})_t&=2\lambda^2\mu^2\beta_i^2e^{2\mu(\beta+6m)}\gamma\gamma_t+\lambda\mu^2\beta_i^2e^{\mu(\beta+6m)}\gamma_t+\lambda\mu\beta_{ii}e^{\mu(\beta+6m)}\gamma_t, \\ \label{eq:iden_li_lii_t}
&=: M_i
\end{align}

On the other hand, after a long but straightforward computation, we get from \eqref{eq:deriv_weights} that
\begin{equation}\label{eq:def_Ali_Alii}
A_i\ell_i+A\ell_{ii}=P^{(1)}_i+P^{(2)}_i
\end{equation}
where 
\begin{align}
P_i^{(1)}&:=-\lambda^2\alpha\gamma_t\gamma\mu^2\beta_i^2e^{\mu(\beta+6m)}-\lambda^2\alpha\gamma_t\gamma\mu\beta_{ii}e^{\mu(\beta+6m)}-\lambda^2\gamma_t\gamma\mu^2\beta_i^2e^{2\mu(\beta+6m)}, \label{eq:def_Pi_1} \\ \notag
P_i^{(2)}&:=\sum_{k}\left[3\lambda^3\mu^4\xi^3\beta_k^2\beta_i^2+2\lambda^2\mu^4\xi^2\beta_k^2\beta_i^2+\lambda^2\mu^3\xi^2\beta_{kk}\beta_i^2+\lambda^3\mu^3\xi^3\beta_k\beta_i\beta_{ki}\right. \\ \label{eq:def_Pi_2}
&\left.\qquad\qquad + 2\lambda^2\mu^3\xi^2\beta_k\beta_{ki}\beta_i+\lambda^3\mu^3\xi^3\beta_k^2\beta_{ii}+\lambda^2\mu^3\xi^2\beta_k^2\beta_{ii}+\lambda^2\mu^2\xi^2\beta_{kk}\beta_{ii}\right].
\end{align}
In the term $P_i^{(2)}$, we have further simplified the notation by recalling that $\xi=e^{\mu(\beta+6m)}\gamma$. Also observe that we have deliberately put together all the terms containing $\gamma_t$ in the above expression.

 We will use now the term $I\Psi\psi$ in \eqref{eq:F_term} to collect other terms containing $\gamma_t$. Indeed, from the definition of $\Psi$ (see eq. \eqref{eq:definicion_Psi}), we see that this term can be rewritten as
 \begin{align}\notag
 \int_{0}^{T}I\Psi\psi\dt&=\int_{0}^{T}A\Psi\psi^2\dt-2\int_{0}^T\left(\sum_{i,k}\psi_{ii}\ell_{kk}\right)\psi \dt \\ \notag
 &= 2 \int_{0}^{T}\sum_{i}P_i^{(3)}\psi^2\dt -2\int_{0}^{T}\left(\sum_{i,k}(\ell_i^2+\ell_{ii})\ell_{kk}\right)\psi^2\dt \\ \label{eq:iden_I_Psi_psi}
 &\quad -2\int_{0}^{T}\left(\sum_{i,k}\psi_{ii}\ell_{kk}\right)\psi\dt,
 \end{align}
 where
 \begin{equation*}
 P_i^{(3)}:=\lambda^2\alpha\gamma\gamma_t\mu^2\beta_i^2e^{\mu(\beta+6m)}+\lambda^2\alpha\gamma\gamma_t\mu\beta_{ii}e^{\mu(\beta+6m)}.
 \end{equation*}

Hence, from \eqref{eq:defs_F}, \eqref{eq:iden_li_lii_t}, \eqref{eq:def_Ali_Alii}, and \eqref{eq:iden_I_Psi_psi}, we get
\begin{align}\notag
\int_{0}^{T}&(F_2+F_3)\psi^2\dt+\int_{0}^{T}I\Psi\psi\dt \\ \notag
&= \underbrace{\int_0^{T}\left(\sum_{i}\left(-\frac{1}{2}M_i-P_i^{(1)}+P_i^{(3)}\right)\right) \psi^2\dt}_{=: Q_1}+\underbrace{\int_{0}^{T}\sum_i\left(P_i^{(2)}-2\sum_k(\ell_i^2+\ell_{ii})\ell_{kk}\right)\psi^2\dt}_{=:Q_2} \\ \label{eq:def_Q_term}
&\quad-\underbrace{2\int_{0}^{T}\left(\sum_{i,k}\psi_{ii}\ell_{kk}\right)\psi\dt}_{=: Q_3}.
\end{align}

We shall focus on the term $Q_1$. From the definition of $M_i$, \eqref{eq:def_Pi_1} and using that $\xi=e^{\mu(\beta+6m)}\gamma$ and $\varphi=\alpha\gamma$, we see that 
\begin{align}\notag
-\frac{M_i}{2}-P_i^{(1)}+P_i^{(3)}=&-\frac{\gamma_t}{\gamma}\left(2\lambda^2\mu^2\xi^2\beta_i^2+\frac{1}{2}\lambda\mu^2\xi\beta_i^2+\frac{1}{2}\lambda\mu\xi\beta_{ii}\right) \\ \label{eq:terms_gamma_t}
&-\frac{\gamma_t}{\gamma}\left(\lambda^2(-\varphi)\xi\mu^2\beta_i^2+\lambda^2(-\varphi)\xi\mu\beta_{ii}\right).
\end{align}
From the definition of $\gamma$, it is clear that the above expression vanishes on $(T/4,T/2)$. On $(T/2,T)$, we use the fact that there exists $C>0$ such that $|\gamma_t|\leq C\gamma^2$.
Hence, for all $(t,x)\in (T/2,T)\times \dom$, there exists a constant $C>0$ only depending on $\dom$ and $\dom^\prime$ such that
\begin{equation}\label{eq:est_Q1_T2_T}
\abs{\sum_{i}\left(-\frac{M_i}{2}-P_i^{(1)}+P_i^{(3)}\right)}\leq C\lambda^2\mu^2\left(\xi^2+\xi\varphi\right)\leq C\lambda^2\mu^3\xi^3,
\end{equation}
where we have used that $|\varphi\gamma|\leq \mu\xi^2$.

On $(0,T/4)$, we are going to use the fact that $\gamma_t\leq 0$, $\varphi<0$, and $\gamma\in[1,2]$ to deduce that $Q_1$ has the good sign outside $\dom^\prime$. Indeed, from \eqref{eq:prop_weight_beta}, we can find $\mu_2=\mu_2(\alpha,\|\Delta\psi\|_{\infty})$ such that for all $\mu\geq \mu_2\geq \mu_1\geq 0$
\begin{align}\notag 
\sum_{i}&\left[2\lambda^2\mu^2\xi^3\beta_i^2+\frac{1}{2}\lambda\mu^2\xi\beta_i^2+\frac{1}{2}\lambda\mu_\xi\beta_{ii}\right]\\ \label{eq:est_pos_D_D0}
&+\sum_{i}\left[\lambda^2(-\varphi)\xi\mu^2\beta_i^2+\lambda^2(-\varphi)\xi\mu\beta_{ii}\right]\geq c\lambda^2\mu^2|\varphi|\xi, \quad x\in \dom\setminus\ov{\dom^\prime}.
\end{align}
In this way, in a subsequent step, by \eqref{eq:terms_gamma_t}, \eqref{eq:est_pos_D_D0}, we will obtain from $Q_1$ a positive term in $(0,T)\times \dom$ and a localized term at $\dom^\prime$ in the right-hand side of the inequality.

The conclusion of this sub-step is quite classical. For the term $Q_2$ in \eqref{eq:def_Q_term}, we can readily see that the leading term in \eqref{eq:def_Pi_2} is positive. Hence, from \eqref{eq:deriv_weights} and  straightforward computation, we have
\begin{equation}\label{eq:est_Q2}
Q_2\geq \int_{0}^{T}\lambda^3\mu^4\xi^3|\nabla \beta|^4\psi^2\dt-C\int_{0}^{T}\left(\lambda^2\mu^4\xi^2+\lambda^2\mu^3\xi^3+\lambda^3\mu^3\xi^3\right)\psi^2\dt
\end{equation}
for some constant $C=C(\|\nabla\beta\|_\infty,\|D^2 \beta\|_\infty)>0$. As in the previous case, using \eqref{eq:prop_weight_beta} will yield a positive term, a localized term on the right-hand side. The terms with lower powers of $\mu$ and $\lambda$ will be absorbed later. 

Finally, for analyzing $Q_3$, we will use that $-\psi_{ii}\ell_{kk}\psi=-(\psi_i\ell_{kk}\psi)_i+\psi_{i}\ell_{kki}\psi+\psi_i^2\ell_{kk}$. Thus,
\begin{equation}\label{eq:Q3_iden}
Q_3=2\int_{0}^{T}\sum_{i,k}\psi_i^2\ell_{kk}\dt+2\int_{0}^{T}\sum_{i,k}\psi_i\psi\ell_{kki}\dt-2\int_{0}^{T}\sum_{i,k}(\psi_i\ell_{kk}\psi)_i\dt.
\end{equation}
We will leave this term as it is. In the next sub-step we will use it for producing a positive term depending on $|\nabla\psi|$. 

\textit{-Estimates on the gradient of $\psi$}. The last positive term we shall obtain in this step, comes from the second term in \eqref{eq:comb_terms} and the first term in \eqref{eq:Q3_iden}. Using \eqref{eq:deriv_weights} it can be readily seen that
\begin{equation}\label{eq:est_grad_final}
\int_{0}^{T}\sum_{i,j}\left(2\psi_i^2\ell_{jj}+c^{ij}\psi_i\psi_j\right)\dt \geq \int_{0}^{T} \lambda\mu^2\xi|\nabla\beta|^2|\nabla\psi|^2 - C\int_{0}^{T}\lambda\mu\xi|\nabla \psi|^2,
\end{equation}
for some $C>0$ only depending on $\dom$ and $\dom^\prime$. From here, using the properties of $\beta$ we will obtain a positive term and a localized term in $\dom^\prime$.

From the second term in \eqref{eq:Q3_iden} and the fact that 
\begin{equation*}
\ell_{kki}=2\lambda\xi\mu^2\beta_k\beta_{ki}+\lambda\xi\mu^3\beta_{k}^2\beta_i+\lambda\xi\mu\beta_{kki}+\lambda\xi\mu^2\beta_{kk}\beta_i,
\end{equation*}
we can use Cauchy-Schwarz and Young inequalities to deduce
\begin{equation}\label{eq:est_cross_grad_sol}
\int_{0}^{T}\sum_{i,k}\psi_i\psi\ell_{kki}\dt \geq - C\int_{0}^{T}\mu^2|\nabla\psi|^2\dt-C\int_0^{T}\lambda^2\mu^4\xi^2|\psi|^2\dt.
\end{equation}
Notice that the term containing $\nabla\psi$ does not have any power for $\lambda$ so it can be absorbed later. 

The last term in \eqref{eq:Q3_iden} is left as it is, since by the divergence theorem we will see later that this term is actually 0.

\textbf{Step 3. Towards the Carleman estimate.} We begin by integrating \eqref{eq:point_iden_lhs}--\eqref{eq:point_last_rhs} in $\dom$ and take expectation in both sides of the identity. Taking into account the estimates obtained in the previous step, i.e., \eqref{eq:est_final_t0}, \eqref{eq:est_F1},  \eqref{eq:def_Q_term}, and \eqref{eq:est_Q1_T2_T}--\eqref{eq:est_cross_grad_sol}, we get
\begin{align}\notag 
c_1&\esp\left(\int_{\dom}|\nabla\psi(0)|^2\dx+\int_{\dom}\lambda^2\mu^3e^{2\mu(6m+1)}|\psi(0)|^2\dx\right)+c\esp\left(\int_0^{T/4}\!\!\!\!\int_{\dom\setminus{\dom^\prime}}\lambda^2\mu^2\xi |\varphi||\gamma_t| |\psi|^2\dx\dt\right) \\ \notag
&\quad +\esp\left(\int_{Q_T}\lambda^3\mu^4\xi^3|\nabla\beta|^4|\psi|^2\dx\dt+\int_{Q_T}\lambda\mu^2\xi|\nabla\beta|^2|\nabla\psi|^2\dx\dt\right)\\ \notag 
&\quad +\esp\left(\int_{Q_T}I^2\dx\dt\right) +\frac{1}{2}\esp\left(\int_{Q_T}\sum_{i}(\d\psi_i)^2\dx\right) \\ \label{eq:first_car_estimate}
&\leq \esp\left(\int_{Q_T}\theta I(\d{z}+\Delta z \dt)\dx\right)+\frac{1}{2}\esp\left(\int_{Q_T} A(\d{\psi})^2\dx\right) + \mathcal{BT}+\mathcal{R},
\end{align}
where 
\begin{align}\label{eq:boundary_terms}
\mathcal{BT}&:=2\esp\left(\int_{Q_T}\sum_{i,k}(\psi_i\ell_k\psi)_i\dx\dt\right)-\esp\left(\int_{Q_T}\sum_{i}(\psi_id\psi)_i\right)-\esp\left(\int_{Q_T}\nabla\cdot V\dx\dt \right), \\ \label{eq:residual}
\mathcal R&:=C\esp\left(\int_{Q_T}\left[\lambda^2\mu^3\xi^3+\lambda^2\mu^4\xi^2+\lambda^3\mu^3\xi^3\right]|\psi|^2\dx\dt+\int_{Q_T}\left[\mu^2+\lambda\mu\xi\right]|\nabla\psi|^2\dx\dt\right).
\end{align}
We remark that the positive constants $c_1$ and $C$ in \eqref{eq:first_car_estimate}--\eqref{eq:residual} only depend on $\dom$ and $\dom^\prime$, while $c>0$ depends only on $\dom,\dom^\prime$ and $\alpha$ (see \eqref{eq:prop_weight_beta}). 

We proceed to estimate the rest of the terms. We begin with those gathered on $\mathcal{BT}$, defined in \eqref{eq:boundary_terms}.												 It is clear that $z=0$ on $\Sigma_T$ implies $\psi=0$ on $\Sigma_T$. Moreover, $\psi_i=\frac{\partial \psi}{\partial \nu}\nu^{i}$, with $\nu=(\nu_1,\ldots,\nu_N)$ being the unit outward normal vector of $\dom$ at $x\in\partial \dom$. Also, by the construction of the weight $\beta$, we have
\begin{equation*}
\ell_i=\lambda\mu\xi\psi_i=\lambda\mu\xi\frac{\partial \psi}{\partial \nu}\nu^{i} \quad\text{and}\quad \frac{\partial \psi}{\partial \nu}<0, \quad\text{on } \Sigma_T.
\end{equation*}
Hence, it is not difficult to see that using divergence theorem we have
\begin{gather*}
2\esp\left(\int_{Q_T}(\psi_i\ell_k\psi)_i\dx\dt\right)=2\esp\left(\int_{\Sigma_T}\sum_{i,k}\psi_i\ell_k\psi \nu^{i}\dx\dt\right)=0, \\
-\esp\left(\int_{Q_T}\sum_i(\psi_i\d{\psi})_i\right)=-\esp\left(\int_{\Sigma_T}\sum_{i}\psi_i\nu_i\d\psi\dx\right)=0,
\end{gather*}
and
\begin{align*}
-\esp\left(\int_{Q_T}\nabla\cdot V\dx\dt\right)&=\esp\left(\int_{\Sigma_T}\sum_{i,j}\left[(2\ell_i\psi_i\psi_j-\ell_i\psi_j\psi_j)+A\ell_i\psi^2\right]\nu^j\dx\dt\right) \\
&= \esp\left(\int_{\Sigma_T}\lambda\mu\xi\frac{\partial \beta}{\partial\nu}\left(\frac{\partial z}{\partial \nu}\right)^2\sum_{i,j}\left(\nu^{i}\nu^j\right)^2\dx\dt\right)\leq 0.
\end{align*}
Thus, we get
\begin{equation}\label{eq:est_boundary}
\mathcal{BT}\leq 0.
\end{equation}

For the following three terms, we will use the change of variables $\psi=\theta z$ and the fact that $z$ solves system \eqref{eq:system_z}. First, we see that
\begin{equation}\label{eq:est_positive_ov_zi}
\esp\left(\int_{Q_T}\sum_i(\d\psi_i)^2\dx\right)=\esp\left(\int_{Q_T}\theta^2\sum_{i}\left(\ov{z}_i+\ell_i\ov{z}\right)^2\dx\dt\right)\geq 0.
\end{equation}
In the same spirit, using the equation verified by $z$ and Cauchy-Schwarz and Young inequalities, we get
\begin{equation}\label{eq:est_rhs_F}
\esp\left(\int_{Q_T}\theta I(\d{z}+\Delta z \dt)\dx\right)\leq \frac{1}{2}\esp\left(\int_{Q_T}I^2\dx\dt\right)+\frac{1}{2}\esp\left(\int_{Q_T}\theta^2|\Xi|^2\dx\dt\right).
\end{equation}
Lastly, from \eqref{eq:deriv_weights} and the fact that $|\varphi_t|\leq C\lambda\mu\xi^3$ for $(t,x)\in(0,T)\times\dom$, a direct computation shows that
\begin{equation}\label{eq:est_A_ov_z}
\esp\left(\int_{Q_T}A(\d\psi)^2\dx\right)=\esp\left(\int_{Q_T}\theta^2A|\ov{z}|^2\dx\dt\right)\leq C\left(\int_{Q_T}\theta^2\lambda^2\mu^2\xi^3|\ov{z}|^2\dx\dt\right).
\end{equation}

Using that $\inf_{x\in\dom\setminus\ov{\dom^\prime}}|\nabla\beta|\geq \alpha>0$, we can combine estimate \eqref{eq:first_car_estimate} with \eqref{eq:est_boundary}--\eqref{eq:est_A_ov_z} to deduce
\begin{align*}\notag 
&\esp\left(\int_{\dom}|\nabla\psi(0)|^2\dx+\int_{\dom}\lambda^2\mu^3e^{2\mu(6m+1)}|\psi(0)|^2\dx\right)+\esp\left(\int_0^{T/4}\!\!\!\!\int_{\dom}\lambda^2\mu^2\xi |\varphi||\gamma_t| |\psi|^2\dx\dt\right) \\ \notag
&\quad +\esp\left(\int_{Q_T}\lambda^3\mu^4\xi^3|\psi|^2\dx\dt+\int_{Q_T}\lambda\mu^2\xi|\nabla\psi|^2\dx\dt\right) +\frac{1}{2}\esp\left(\int_{Q_T}I^2\dx\dt\right)  \\ \notag
&\leq C\esp\left(\int_0^{T/4}\!\!\!\!\int_{\dom^\prime}\lambda^2\mu^2\xi |\varphi||\gamma_t| |\psi|^2\dx\dt+\int_{0}^{T}\!\!\!\int_{\dom^\prime}\lambda^3\mu^4\xi^3|\psi|^2\dx\dt+\int_{0}^{T}\!\!\!\int_{\dom^\prime}\lambda\mu^2\xi|\nabla\psi|^2\dx\dt\right) \\
&\quad +  C\mathcal R + C \esp\left(\int_{Q_T}\theta^2|\Xi|^2\dx\dt+\int_{Q_T}\theta^2\lambda^2\mu^2\xi^3|\ov{z}|^2\dx\dt\right),
\end{align*}
for some $C>0$ only depending on $\dom$, $\dom^\prime$ and $\alpha$. We observe that, unlike the traditional Carleman estimate with weight vanishing at $t=0$ and $t=T$, we have three local integrals, one of those being only for $t\in(0,T/4)$. We will handle this in the following step. 

Also notice that all of the terms in $\mathcal R$ have lower powers of $\lambda$ and $\mu$, thus, we immediately see that there exists some $\mu_3\geq \mu_2$ and $\lambda_1\geq C$ such that, for all $\mu\geq \mu_3$ and $\lambda\geq \lambda_1$ 
\begin{align}\notag 
&\esp\left(\int_{\dom}|\nabla\psi(0)|^2\dx+\int_{\dom}\lambda^2\mu^3e^{2\mu(6m+1)}|\psi(0)|^2\dx\right)+\esp\left(\int_0^{T/4}\!\!\!\!\int_{\dom}\lambda^2\mu^2\xi |\varphi||\gamma_t| |\psi|^2\dx\dt\right) \\ \notag
&\quad +\esp\left(\int_{Q_T}\lambda^3\mu^4\xi^3|\psi|^2\dx\dt+\int_{Q_T}\lambda\mu^2\xi|\nabla\psi|^2\dx\dt\right)   \\ \notag
&\leq C\esp\left(\int_0^{T/4}\!\!\!\!\int_{\dom^\prime}\lambda^2\mu^2\xi |\varphi||\gamma_t| |\psi|^2\dx\dt+\int_{0}^{T}\!\!\!\int_{\dom^\prime}\lambda^3\mu^4\xi^3|\psi|^2\dx\dt+\int_{0}^{T}\!\!\!\int_{\dom^\prime}\lambda\mu^2\xi|\nabla\psi|^2\dx\dt\right) \\ 
&\quad + C \esp\left(\int_{Q_T}\theta^2|\Xi|^2\dx\dt+\int_{Q_T}\theta^2\lambda^2\mu^2\xi^3|\ov{z}|^2\dx\dt\right).
\label{eq:est_almost0}
\end{align}

\smallskip
\textbf{Step 4. Last arrangements and conclusion.} As usual, the last steps in Carleman strategies consist in removing the local term containing the gradient of the solution and coming back to the original variable. We will see that the original strategy also helps to remove the local term in $(0,T/4)$. 

First, using that $z_i=\theta^{-1}(\psi_i-\ell_i\psi)$, it is not difficult to see that $\theta^2|\nabla z|^2\leq 2|\nabla \psi|^2+2C\lambda^2\mu^2\xi^2|\psi|^2$ for some $C>0$ only depending on $\dom$ and $\dom^\prime$, hence from \eqref{eq:est_almost0} we have
\begin{align}\notag 
&\esp\left(\int_{\dom}\theta^2(0)|\nabla z(0)|^2\dx+\int_{\dom}\lambda^2\mu^3e^{2\mu(6m+1)}\theta^2(0)|z(0)|^2\dx\right)+\esp\left(\int_0^{T/4}\!\!\!\!\int_{\dom}\theta^2\lambda^2\mu^2\xi |\varphi||\gamma_t| |z|^2\dx\dt\right) \\ \notag
&\quad +\esp\left(\int_{Q_T}\theta^2\lambda^3\mu^4\xi^3|z|^2\dx\dt+\int_{Q_T}\theta^2\lambda\mu^2\xi|\nabla z|^2\dx\dt\right)   \\ \notag
&\leq C\esp\left(\int_0^{T/4}\!\!\!\!\int_{\dom^\prime}\theta^2\lambda^2\mu^2\xi |\varphi||\gamma_t| |z|^2\dx\dt+\int_{0}^{T}\!\!\!\int_{\dom^\prime}\theta^2\lambda^3\mu^4\xi^3|z|^2\dx\dt+\int_{0}^{T}\!\!\!\int_{\dom^\prime}\theta^2\lambda\mu^2\xi|\nabla z|^2\dx\dt\right) \\ \label{eq:est_almost}
&\quad + C \esp\left(\int_{Q_T}\theta^2|\Xi|^2\dx\dt+\int_{Q_T}\theta^2\lambda^2\mu^2\xi^3|\ov{z}|^2\dx\dt\right),
\end{align}
for all $\lambda\geq \lambda_1$ and $\mu\geq \mu_3$. 

We choose a cut-off function $\eta\in C_c^\infty(\dom)$ such that
\begin{equation}\label{eq:cut_off}
0\leq \eta\leq 1, \quad \eta\equiv 1 \text{ in } \dom^\prime, \quad \eta\equiv 0 \text{ in }\dom\setminus\dom_0
\end{equation}
with the additional characteristic that
\begin{equation}\label{eq:prop_nabla_cut}
\frac{\nabla \eta}{\eta^{1/2}}\in L^\infty(\dom)^N.
\end{equation}
This condition can be obtained by taking some $\eta_0\in C_c^\infty(\dom)$ satistying \eqref{eq:cut_off} and defining $\eta=\eta_0^4$. Then $\eta$ will satisfy both \eqref{eq:cut_off} and \eqref{eq:prop_nabla_cut}.

Using It\^{o}'s formula, we compute $\d\left(\theta^2\xi z^2\right)=(\theta^2\xi)_tz^2+2\theta^2\xi z\d{z}+\theta^2\xi(\d{z})^2$
and thus, using the equation verified by $z$, we get
\begin{align}\notag 
\esp&\left(\int_{\dom_0}\theta^2(0)\xi(0)|z(0)|^2\eta \dx\right)+2\esp\left(\int_0^T\!\!\!\int_{\dom_0}\theta\theta_t\xi |z|^2\eta\dx\dt\right)+2\esp\left(\int_{0}^{T}\!\!\!\int_{\dom_0}\theta^2\xi|\nabla z|^2\eta\dx\dt\right)\\ \notag
&+\esp\left(\int_{0}^{T}\!\!\!\int_{\dom_0}\theta^2\xi|\ov{z}|^2\eta\dx\dt\right) =-\esp\left(\int_0^T\!\!\!\int_{\dom_0}\theta^2\xi_t |z|^2\eta \dx\dt\right)-2\esp\left(\int_0^T\!\!\!\int_{\dom_0}\theta^2\xi z \Xi \eta \dx\dt \right) \\ \label{eq:iden_local}
&-2\esp\left(\int_{0}^{T}\!\!\!\int_{\dom_0}\theta^2\xi\nabla\eta\cdot\nabla z z \dx\dt\right)-2\esp\left(\int_{0}^{T}\!\!\!\int_{\dom_0}\nabla(\theta^2\xi)\cdot\nabla z z \eta\dx\dt\right).
\end{align}

We readily see that the first and last terms in the left-hand side of \eqref{eq:iden_local} are positive, so they can be dropped. Also, notice that using the properties of $\eta$, the third term gives (up to the constants $\mu$ and $\lambda$) the local term containing $|\nabla z|$. 

We shall focus on the second term on the left-hand side of \eqref{eq:iden_local}. Similar to Step 2 above, we analyze it on different time intervals. Obviously, for $t\in(T/4,T/2)$ this term vanishes since $\gamma_t=0$. For $t\in(0,T/4)$, we notice that $\theta\theta_t=\theta^2\lambda\varphi\frac{\gamma_t}{\gamma}\xi$ and since $\gamma_t\leq 0$, $\varphi<0$ and $\gamma\in[1,2]$, this yields a positive term. Lastly, in the interval $(T/2,T)$, we use that $|\varphi_t|\leq C\lambda\mu \xi^3$ to obtain the bound $|\theta_t|\leq C\theta\lambda^2\mu\xi^3$. Summarizing, we have
\begin{align}\notag 
2\esp&\left(\int_{0}^{T}\!\!\!\int_{\dom_0}\theta\theta_t\xi|z|^2\eta\dx\dt\right)\\ \label{eq:est_positive_t}
&\geq \esp\left(\int_{0}^{T/4}\!\!\!\int_{\dom_0}\theta^2\lambda\xi|\gamma_t||\varphi||z|^2\eta\dx\dt\right)- C\esp\left(\int_{T/2}^{T}\int_{\dom_0}\theta^2\lambda^2\mu\xi^3|z|^2\eta\dx\dt\right).
\end{align}

Let us estimate each term on the right-hand side of \eqref{eq:iden_local}. For the first one, using that $|\xi_t|\leq C\lambda\mu\xi^3$ for all $(t,x)\in(0,T)\times \dom$, we get
\begin{equation}\label{eq:est_xi_t}
\abs{\esp\left(\int_0^T\!\!\!\int_{\dom_0}\theta^2\xi_t |z|^2\eta \dx\dt\right)}\leq C\esp\left(\int_{0}^{T}\!\!\!\int_{\dom_0}\theta^2\lambda\mu\xi^3|z|^2\eta \dx\dt\right).
\end{equation}
For the second one, using Cauchy-Schwarz and Young inequalities yields
\begin{align}\notag 
&\abs{\esp\left(\int_0^T\!\!\!\int_{\dom_0}\theta^2\xi z \Xi \eta \dx\dt \right)}\\ \label{eq:est_cross_Fz}
&\quad \leq \frac{1}{2}\esp\left(\int_{0}^{T}\!\!\!\int_{\dom_0}\theta^2\lambda^{-1}\mu^{-2}|\Xi|^2\eta\dx\dt\right)+\frac{1}{2}\esp\left(\int_{0}^{T}\!\!\!\int_{\dom_0}\theta^2\lambda\mu^2\xi^2|z|^2\eta\dx\dt\right).
\end{align}
For the third one, we will use property \eqref{eq:prop_nabla_cut} and Cauchy-Schwarz and Young inequalities to deduce that
\begin{align}\notag 
&\abs{\esp\left(\int_{0}^{T}\!\!\!\int_{\dom_0}\theta^2\xi\nabla\eta\cdot\nabla z z \dx\dt\right)}\\ \label{eq:est_nablaz_z}
&\qquad \leq \epsilon \esp\left(\int_{0}^{T}\!\!\!\int_{\dom_0}\theta^2\xi |\nabla z|^2\eta\dx\dt\right)+C(\epsilon)\esp\left(\int_{0}^{T}\!\!\!\int_{\dom_0}\theta^2\xi|z|^2\dx\dt\right)
\end{align}
for any $\epsilon>0$. For the last term, using that $|\nabla(\theta^2\xi)|\leq C\theta^2\lambda\mu\xi^2$ and arguing as above, we get
\begin{align}\notag 
&\abs{\esp\left(\int_{0}^{T}\!\!\!\int_{\dom_0}\nabla(\theta^2\xi)\cdot\nabla z z \eta\dx\dt\right)}\\ \label{eq:est_nablaweight_nablaz}
&\qquad \leq \epsilon \esp\left(\int_{0}^{T}\!\!\!\int_{\dom_0}\theta^2\xi |\nabla z|^2\eta\dx\dt\right) + C(\epsilon)\esp\left(\int_{0}^{T}\!\!\!\int_{\dom_0}\theta^2\xi^3\mu^2\lambda^2|z|^2\eta\dx\dt\right).
\end{align}
Therefore, taking $\epsilon=\frac{1}{2}$ and using estimates \eqref{eq:est_positive_t}--\eqref{eq:est_nablaweight_nablaz} together with the properties of the cut-off $\eta$, we get
\begin{align}\notag 
\esp&\left(\int_{0}^{T/4}\!\!\!\int_{\dom^\prime}\theta^2\lambda\xi|\gamma_t||\varphi||z|^2\dx\dt\right)+\esp\left(\int_{0}^{T}\!\!\!\int_{\dom^\prime}\theta^2\xi|\nabla z|^2\dx\dt\right) \\ \label{eq:est_locales}
&\leq C\esp\left(\int_{0}^{T}\int_{\dom_0}\theta^2(\lambda^2\mu\xi^3+\lambda\mu^2\xi^2+\lambda^2\mu^2\xi^3)|z|^2\dx\dt\right)+C\esp\left(\int_{Q_T}\theta^2\lambda^{-1}\mu^{-2}|\Xi|^2\dx\dt\right).
\end{align}
As usual, we have paid the price of estimating locally the gradient by enlarging a little bit the observation domain. Notice that this procedure gives us the local estimate in $(0,T/4)$ by using the properties of the weight function $\varphi$ and $\gamma_t$. Finally, the desired estimate follows by multiplying both sides of \eqref{eq:est_locales} by $\lambda\mu^2$ and using the result to bound in the right-hand side of \eqref{eq:est_almost}. We conclude the proof by setting $\mu_0=\mu_3$ and $\lambda_0=\lambda_1$. 
\end{proof}

\subsection{A controllability result for a linear forward stochastic heat equation with two source terms and two controls}

In this section, we will prove a controllability result for a linear forward equation. More precisely, recall the equation defined in \eqref{eq:sys_forward_source_intro}
\begin{equation}\label{eq:sys_forward_source}
\begin{cases}
\d{y}=(\Delta y + F+ \chi_{\dom_0}h)\dt + (G+H)\d{W}(t) &\text{in }Q_T, \\
y=0 &\text{on } \Sigma_T, \\
y(0)=y_0 &\text{in }\dom.
\end{cases}
\end{equation}

In \eqref{eq:sys_forward_source}, $(h,H)\in L^2_{\mathcal F}(0,T;L^2(\dom_0))\times L^2_{\mathcal F}(0,T;L^2(\dom))$ is a pair of controls and $F, G$ are given source terms in $L_{\mathcal F}^2(0,T;L^2(\dom))$. Observe that given $y_0\in L^2(\Omega,\mathcal F_0; L^2(\dom))$ and the aforementioned regularity on the controls and source terms, system \eqref{eq:sys_forward_source} admits a unique solution $y\in \mathcal W_T$, see  \cite[Theorem 2.7]{LZ19}.

Under the notation of \Cref{sec:new_carleman}, let us set the parameters $\lambda$ and $\mu$ to a fixed value sufficiently large, such that inequality \eqref{eq:car_backward_0} holds true.
We define the space 
\begin{align}
&\mathcal S_{\lambda,\mu}=\Bigg\{(F,G)\in [L^2_{\mathcal F}(0,T;L^2(\dom))]^2: \notag\\
& \quad  \left[\esp\left(\int_{Q_T}\theta^{-2}\lambda^{-3}\mu^{-4}\xi^{-3}|F|^2\dx\dt\right)+\esp\left(\int_{Q_T}\theta^{-2}\lambda^{-2}\mu^{-2}\xi^{-3}|G|^2\dx\dt \right)\right]^{1/2}<+\infty \Bigg\},\label{eq:defslambdamu}
\end{align}
endowed with the canonical norm. 

Our linear controllability result reads as follows.
\begin{theo}\label{teo:contr_forward_source}
For any initial datum $y_0\in L^2(\Omega,\mathcal F_0;L^2(\dom))$ and any source terms $(F,G)\in \mathcal S_{\lambda,\mu}$, there exists a pair of controls $(\widehat{h},\widehat{H})\in L^2_{\mathcal F}(0,T;L^2(\dom_0))\times L^2_{\mathcal F}(0,T;L^2(\dom))$ such that the associated solution $\widehat{y} \in \mathcal{W}_T$ to system \eqref{eq:sys_forward_source} satisfies $\widehat{y}(T)=0$ in $\dom$, a.s. Moreover, the following estimate holds
\begin{align}\notag 
\esp&\left(\int_{Q_T}\theta^{-2}|\widehat{y}|^2\dx\dt\right)+\esp\left(\int_{0}^{T}\!\!\!\int_{\dom_0}\theta^{-2}\lambda^{-3}\mu^{-4}\xi^{-3}|\widehat{h}|^2\dx\dt\right)+\esp\left(\int_{Q_T}\theta^{-2}\lambda^{-2}\mu^{-2}\xi^{-3}|\widehat{H}|^2\dx\dt\right) \\ \label{eq:est_weighted_spaces_forward}
& \leq C_1\esp\left(\|y_0\|^2_{L^2(\dom)}\right)+C\norme{(F,G)}_{\mathcal S_{\lambda,\mu}}^2,
\end{align}
where $C_1>0$ is a constant depending on $\dom,\dom_0,\lambda,\mu$ and $C>0$ only depends on $\dom$ and $\dom_0$.
\end{theo}
\begin{rmk}
\label{rmk:uniquetrajectory}
Using classical arguments, see for instance \cite[Proposition 2.9]{LLT13}, from \Cref{teo:contr_forward_source} one can construct  a linear continuous mapping that associates every initial datum $y_0\in L^2(\Omega,\mathcal F_0;L^2(\dom))$ and every source terms $(F,G) \in \mathcal{S}_{\lambda, \mu}$, to a trajectory $(\widehat{y}, \widehat{h}, \widehat{H})$ such that $\widehat{y}(T)=0$ in $\dom$, a.s. and \eqref{eq:est_weighted_spaces_forward} holds.
\end{rmk}

The proof of \Cref{teo:contr_forward_source} is based on a classical duality method, called penalized Hilbert Uniqueness Method, which ideas can be traced back to the seminal work \cite{GL94}. The general strategy consists in three steps:
\begin{itemize}
\item[-] \textit{Step 1.} Construct a family of optimal approximate-null control problems for system
\eqref{eq:sys_forward_source}.
\item[-] \textit{Step 2.} Obtain a uniform estimate for the approximate solutions in terms of the data of the problem, i.e., the initial datum $y_0$ and the source terms $F$ and $G$.
\item[-] \textit{Step 3.} A limit process to derive the desired null-controllability result. 
\end{itemize}

We shall mention that in the stochastic setting, similar strategies have been used for deducing controllability results and Carleman estimates for forward and backward equations, see e.g. \cite{liu14,Yan18,LY19}.

In what follows, $C$ will denote a generic positive constant possibly depending on $\dom,\dom_0$, but never on the parameters $\lambda$ and $\mu$.

\begin{proof}[Proof of \Cref{teo:contr_forward_source}] We follow the steps described above. 

\smallskip
\textbf{Step 1}. For any $\epsilon>0$, let us consider the weight function $\theta_\epsilon(t)$ given by
\begin{equation*}
\gamma_\epsilon(t):=
\begin{cases}
\gamma_\epsilon(t)=1+(1+\frac{4t}{T})^\sigma, \quad t\in[0,T/4], \\
\gamma_\epsilon(t)=1, \quad t\in [T/4,T/2+\epsilon], \\
\gamma_\epsilon(t)=\gamma(t-\epsilon), \quad t\in [T/2+\epsilon,T], \\
\sigma \textnormal{ as in \eqref{def:sigma}}. 
\end{cases}
\end{equation*}
Defined in this way, it is not difficult to see that $\gamma$ does not blow up as $t\to T^{-}$ and that $\gamma_\epsilon(t)\leq \gamma(t)$ for $t\in[0,T]$. With this new function, we set the weight $\varphi_\epsilon$ as in \eqref{eq:weights_0} by replacing the function $\gamma$ by $\gamma_\epsilon$. In the same manner, we write $\theta_\epsilon=e^{\lambda\varphi_\epsilon}$.

With this notation, we introduce the functional
\begin{align}\notag
J_\epsilon(h,H):=&\frac{1}{2}\esp\left(\int_{Q_T}\theta_\epsilon^{-2}|y|^2\dx\dt\right)+\frac{1}{2}\esp\left(\int_{0}^T\intzero \theta^{-2}\lambda^{-3}\mu^{-4}\xi^{-3}|h|^2\dx\dt\right)\\ \label{eq:func}
&+\frac{1}{2}\esp\left(\int_{Q_T}\theta^{-2}\lambda^{-2}\mu^{-2}\xi^{-3}|H|^2\dx\dt\right)+\frac{1}{2\epsilon}\esp\left(\int_{\dom}|y(T)|^2\dx\right)
\end{align}
and consider the minimization problem
\begin{equation}\label{eq:prob_min}
\begin{cases}
\min_{(h,H)\in \mathcal H} J_\epsilon(h,H), \\
\textnormal{subject to equation \eqref{eq:sys_forward_source},} 
\end{cases}
\end{equation}
where \begin{align*}
\mathcal H=&\left\{(h,H)\in L^2_{\mathcal F}(0,T;L^2(\dom)): \right. \\ &\quad \left.
\esp\left(\int_{0}^{T}\!\!\!\int_{\dom_0}\theta^{-2}\lambda^{-3}\mu^{-4}\xi^{-3}|h|^2\dx\dt\right)<+\infty, \ \esp\left(\int_{Q_T}\theta^{-2}\lambda^{-2}\mu^{-2}\xi^{-3}|H|^2\dx\dt\right)<+\infty\right\}.
\end{align*}

It can be readily seen that the functional $J_\epsilon$ is continuous, strictly convex and coercive. Therefore, the minimization problem \eqref{eq:prob_min} admits a unique optimal pair solution that we denote by $(h_\epsilon,H_\epsilon)$. From classical arguments: the Euler-Lagrange equation for \eqref{eq:func} at the minimum $(h_\epsilon,H_\epsilon)$ and a duality argument (see, for instance, \cite{Lio71}), the pair $(h_\epsilon,H_\epsilon)$ can be characterized as 
\begin{equation}\label{eq:h_epsilon}
h_\epsilon=-\chi_{\dom_0}\theta^2\lambda^3\mu^4\xi^3 z_\epsilon, \quad H_\epsilon=-\theta^2\lambda^2\mu^2\xi^3Z_{\epsilon} \quad\text{in Q}, \quad a.s.,
\end{equation}
where the pair $(z_\epsilon,Z_\epsilon)$ verifies the backward stochastic equation
\begin{equation}\label{eq:q_optim}
\begin{cases}
\d{z_\epsilon}=(-\Delta z_\epsilon-\theta_\epsilon^{-2}y_\epsilon)\dt+ Z_\epsilon \d{W}(t)&\text{in }Q_T, \\
z_\epsilon=0 &\text{on }\Sigma_T, \\
z_\epsilon(T)=\frac{1}{\epsilon}y_\epsilon(T) &\text{in }\dom,
\end{cases}
\end{equation}
and where $(y_\epsilon,y_\epsilon(0))$ can be extracted from $y_\epsilon$ the solution to \eqref{eq:sys_forward_source} with controls $h=h_\epsilon$ and $H=H_\epsilon$. Observe that since $y_\epsilon\in L^2_{\mathcal F}(\Omega; C([0,T];L^2(\dom)))$ the evaluation of $y_\epsilon$ at $t=T$ is meaningful and \eqref{eq:q_optim} is well-posed for any $\epsilon>0$.

\smallskip
\textbf{Step 2.} Using It\^{o}'s formula, we can compute $\d(y_\epsilon z_\epsilon)$ and deduce
\begin{align*}
\esp\left(\int_{\dom}y_\epsilon(T)z_\epsilon(T)\dx\right)&=\esp\left(\int_{\dom}y_\epsilon(0)z_\epsilon(0)\dx\right)+\esp\left(\int_{Q_T}(\Delta y_\epsilon+F+\chi_{\dom_0}h_\epsilon)z_\epsilon\dx\dt\right) \\
&\quad + \esp\left(\int_{Q_T}(-\Delta z_\epsilon-\theta_\epsilon^{-2}y_\epsilon)y_\epsilon\dx\dt\right)+\esp\left(\int_{Q_T}(H_\epsilon+G)Z_\epsilon \dx\dt \right)
\end{align*}
whence, replacing the initial data of systems \eqref{eq:sys_forward_source}, \eqref{eq:q_optim} and using identity \eqref{eq:h_epsilon}, we get
\begin{align}
\esp&\left(\int_{0}^{T}\!\!\!\int_{\dom_0}\theta^2\lambda^3\mu^4\xi^3|z_\epsilon|^2\dx\dt\right)+\esp\left(\int_{Q_T}\theta^2\lambda^2\mu^2\xi^3|Z_\epsilon|^2\dx\dt\right)\notag\\
&\quad +\esp\left(\int_{Q_T}\theta_\epsilon^{-2}|y_\epsilon|^2\dx\dt\right)+\frac{1}{\epsilon}\esp\left(\int_{\dom}|y_\epsilon(T)|^2\dx\right) \notag \\ \label{eq:id_prod_yq}
&=\esp\left(\int_{\dom}y_0z_\epsilon(0)\dx\right)+\esp\left(\int_{Q_T}F z_\epsilon \dx\dt \right)+\esp\left(\int_{Q_T}G Z_\epsilon \dx\dt \right).
\end{align}

Now, we will use the Carleman estimate in \Cref{thm:carleman_backward_0}. We will apply it to equation \eqref{eq:q_optim} with $\Xi=-\theta^{-2}y_\epsilon$ and $\ov{z}=Z_\epsilon$. Then, after removing some unnecessary terms, we get for any $\lambda$ and $\mu$ large enough
\begin{align}\notag
\esp&\left(\int_{\dom}\lambda^2\mu^3\theta^2(0)|z_\epsilon(0)|^2\dx\right)+\esp\left(\int_{Q_T}\lambda^3\mu^4\xi^3\theta^2| z_\epsilon|^2\dx\dt\right)+\esp\left(\int_{Q_T}\lambda^2\mu^2\xi^3\theta^2| Z_\epsilon |^2\dx\dt\right) \\  \label{eq:car_random}
&\leq C\esp\left(\int_{0}^{T}\!\!\!\int_{\dom_0}\lambda^3\mu^4\xi^3\theta^2|z_\epsilon|^2\dx\dt+\int_{Q_T}\theta^2|\theta_\epsilon^{-2} y_\epsilon|^2\dx\dt+\int_{Q_T}\lambda^2\mu^2\xi^3\theta^2|Z_\epsilon|^2\dx\dt\right).
\end{align}
Notice that we have added an integral of $Z_\epsilon$ on the left-hand side of the inequality. This increases a little bit the constant $C$ on the right-hand side but it is still uniform with respect to $\lambda$ and $\mu$.

In view of \eqref{eq:car_random}, we use Cauchy-Schwarz and Young inequalities in the right-hand side of \eqref{eq:id_prod_yq} to obtain
\begin{align}\notag
\esp&\left(\int_{0}^{T}\!\!\!\int_{\dom_0}\theta^2\lambda^3\mu^4\xi^3|z_\epsilon|^2\dx\dt\right)+\esp\left(\int_{Q_T}\theta^2\lambda^3\mu^4\xi^3| z_\epsilon|^2\dx\dt\right) \\ \notag
&\quad +\esp\left(\int_{Q_T}\theta^{-2}|y_\epsilon|^2\dx\dt\right)+\frac{1}{\epsilon}\esp\left(\int_{\dom}|y_\epsilon(T)|^2\dx\right) \\ \notag
&\leq \delta\left[\esp\left(\int_{\dom}\theta^2(0)\lambda^2\mu^3|z_\epsilon(0)|^2\dx\right)+\esp\left(\int_{Q_T}\theta^2\lambda^3\mu^4\xi^3|z_\epsilon|^2\dx\dt+\int_{Q_T}\theta^2\lambda^2\mu^2\xi^3|Z_\epsilon|^2\dx\dt\right)\right]\\ \notag 
&\quad + C_\delta\left[\esp\left(\int_{\dom}\theta^{-2}(0)\lambda^{-2}\mu^{-3}|y_0|^2\dx\right)+\esp\left(\int_{Q_{T}}\theta^{-2}\lambda^{-3}\mu^{-4}\xi^{-3}|F|^2\dx\dt\right) \right. \\ \label{eq:des_prod_yq}
&\qquad \qquad \left. + \esp\left(\int_{Q_{T}}\theta^{-2}\lambda^{-2}\mu^{-2}\xi^{-3}|G|^2\dx\dt\right) \right]
\end{align}
for any $\delta>0$. Using inequality \eqref{eq:car_random} to estimate in the right-hand side of \eqref{eq:des_prod_yq} and the fact that $\theta^2 \theta_\epsilon^{-2}\leq 1$ for all $(t,x)\in Q_T$, we obtain, after taking $\delta>0$ small enough, that
\begin{align*}
\esp&\left(\int_{0}^{T}\!\!\!\int_{\dom_0}\theta^2\lambda^3\mu^4\xi^3|z_\epsilon|^2\dx\dt\right)+\esp\left(\int_{Q_T}\theta^2\lambda^2\mu^2\xi^3|Z_\epsilon|^2\dx\dt\right)\\
&\quad +\esp\left(\int_{Q_T}\theta_\epsilon^{-2}|y_\epsilon|^2\dx\dt\right)+\frac{1}{\epsilon}\esp\left(\int_{\dom}|y_\epsilon(T)|^2\dx\right) \\ 
&\leq C\left[\esp\left(\int_{\dom}\theta^{-2}(0)\lambda^{-2}\mu^{-3}|y_0|^2\dx\right)+\esp\left(\int_{Q_{T}}\theta^{-2}\lambda^{-3}\mu^{-4}\xi^{-3}|F|^2\dx\dt\right) \right. \\
&\qquad\quad  \left. +\esp\left(\int_{Q_{T}}\theta^{-2}\lambda^{-2}\mu^{-2}\xi^{-3}|G|^2\dx\dt\right) \right].
\end{align*}
Recalling the characterization of the optimal control $h_\epsilon$ in \eqref{eq:h_epsilon} we obtain
\begin{align}\notag 
\esp&\left(\int_{0}^{T}\!\!\!\int_{\dom_0}\theta^{-2}\lambda^{-3}\mu^{-4}\xi^{-3}|h_\epsilon|^2\dx\dt\right)+\esp\left(\int_{Q_T}\theta^{-2}\lambda^{-2}\mu^{-2}\xi^{-3}|H_\epsilon|^2\dx\dt\right)\\ \notag
&\quad +\esp\left(\int_{Q_T}\theta_\epsilon^{-2}|y_\epsilon|^2\dx\dt\right)+\frac{1}{\epsilon}\esp\left(\int_{\dom}|y_\epsilon(T)|^2\dx\right) \\ \notag
&\leq C\left[\esp\left(\int_{\dom}\theta^{-2}(0)\lambda^{-2}\mu^{-3}|y_0|^2\dx\right)+\esp\left(\int_{Q_{T}}\theta^{-2}\lambda^{-3}\mu^{-4}\xi^{-3}|F|^2\dx\dt\right) \right. \\  \label{eq:iden_uniform_final}
&\qquad\quad  \left. +\esp\left(\int_{Q_{T}}\theta^{-2}\lambda^{-2}\mu^{-2}\xi^{-3}|G|^2\dx\dt\right) \right].
\end{align}
Observe that the right-hand side of \eqref{eq:iden_uniform_final} is well-defined and finite since $\theta^{-2}(0)<+\infty$ and the source terms $(F,G)$ belongs to $ \mathcal{S}_{\lambda, \mu}$, defined in \eqref{eq:defslambdamu}.

\smallskip
\textbf{Step 3.} Since the right-hand side of \eqref{eq:iden_uniform_final} is uniform with respect to $\epsilon$, we readily deduce that there exists $(\widehat{h},\widehat{y},\widehat{Y})$ such that
\begin{equation}\label{eq:weak_conv}
\begin{cases}
h_\epsilon\weakly \widehat{h} &\textnormal{weakly in } L^2(\Omega\times(0,T);L^2(\dom_0)), \\
H_\epsilon\weakly \widehat{H} &\textnormal{weakly in } L^2(\Omega\times(0,T);L^2(\dom)), \\
y_{\epsilon}\weakly \widehat{y} &\textnormal{weakly in } L^2(\Omega\times(0,T);L^2(\dom)).
\end{cases}
\end{equation}

We claim that $\widehat{y}$ is the solution to \eqref{eq:sys_forward_source} associated to $(\widehat{h},\widehat{H})$. To show this, let us denote by $\tilde{y}$ the unique solution in $L^2_{\mathcal F}(0,T;C([0,T];L^2(\dom)))\cap L^2_{\mathcal F}(0,T;H_0^1(\dom))$ to \eqref{eq:sys_forward_source} with controls $(\widehat{h},\widehat{H})$. For any $m\in L^2_{\mathcal F}(0,T;L^2(\dom))$, we consider $(z,Z)$ the unique solution to the backward equation
\begin{equation}\label{eq:z_gen}
\begin{cases}
\d{z}=(-\Delta z-m)\dt+Z\d{W}(t) &\text{in }Q_T, \\
z=0 &\text{on }\Sigma_T, \\
z(T)=0 &\text{in } \dom.
\end{cases}
\end{equation}
Then, using It\^{o}'s formula, we compute the duality between \eqref{eq:z_gen} and \eqref{eq:sys_forward_source} associated to $(h,H)=(h_\epsilon,H_\epsilon)$ and $(h,H)=(\widehat{h},\widehat{H})$, respectively. We have
\begin{align}\notag 
-\esp\left(\int_{\dom}y_0z(0)\dx\right)&=-\esp\left(\int_{Q_T} m y_\epsilon \dx\dt\right)+\esp\left(\int_{Q_T}F z\dx\dt\right)+\esp\left(\int_{Q_T}G Z\dx\dt\right) \\ \label{eq:duality_eps}
&\quad +\esp\left(\int_{0}^{T}\!\!\!\int_{\dom_0}h_\epsilon z\dx\dt\right) + \esp\left(\int_{Q_T}H_\epsilon z\dx\dt\right)
\end{align}
and
\begin{align}\notag 
-\esp\left(\int_{\dom}y_0z(0)\dx\right)&=-\esp\left(\int_{Q_T} m \tilde y \dx\dt\right)+\esp\left(\int_{Q_T}F z\dx\dt\right)+\esp\left(\int_{Q_T}G Z\dx\dt\right) \\ \label{eq:duality_tilde}
&\quad +\esp\left(\int_{0}^{T}\!\!\!\int_{\dom_0}\widehat{h} z\dx\dt\right) + \esp\left(\int_{Q_T}\widehat{H} z\dx\dt\right).
\end{align}
Then, using \eqref{eq:weak_conv} in \eqref{eq:duality_eps} to pass to the limit $\epsilon\to 0$ and subtracting the result from \eqref{eq:duality_tilde}, we get $\tilde y=\widehat{y}$ in $Q_T$, a.s. 

To conclude, we notice from \eqref{eq:iden_uniform_final} that $\widehat{y}(T)=0$ in $\dom$, a.s. Also, from the weak convergence \eqref{eq:weak_conv}, Fatou's lemma and the uniform estimate \eqref{eq:iden_uniform_final} we deduce \eqref{eq:est_weighted_spaces_forward}. This ends the proof.
\end{proof}

\subsection{Proof the nonlinear result for the forward equation}

Now, we are in position to prove \Cref{th:semilinear_forward}. To this end, let us fix the parameters $\lambda$ and $\mu$ in \Cref{teo:contr_forward_source} to a fixed value sufficiently large. Recall that in turn, these parameters come from \Cref{thm:carleman_backward_0} and should be selected as $\lambda \geq \lambda_0$ and $\mu\geq \mu_0$ for some $\lambda_0\geq 1$ and $\mu_0\geq 1$, so there is no contradiction. 

Note that at this point, we have preserved explicitly the parameters $\lambda$ and $\mu$ in the controllability result of \Cref{teo:contr_forward_source}. This was possible due to the selection of the weight $\theta$ in the Carleman estimate \eqref{eq:car_backward_0}, which allows to have a term depending on $z(0)$ in the left-hand side.

\begin{proof}[Proof of \Cref{th:semilinear_forward}]
Let us consider nonlinearities $f$ and $g$ fulfilling \eqref{eq:f_g_zero} and \eqref{eq:UniformLipschitzf_forw} and define the nonlinear map
\begin{equation*}
\mathcal N: (F,G)\in \mathcal S_{\lambda,\mu}\mapsto (f(\omega,t,x,y),g(\omega,t,x,y))\in \mathcal S_{\lambda,\mu},
\end{equation*}
where $y$ is the trajectory of \eqref{eq:sys_forward_source} associated to the data $y_0$, $F$, and $G$, see \Cref{teo:contr_forward_source} and \Cref{rmk:uniquetrajectory}. In what follows, to abridge the notation, we simply write $f(y)$ and $g(y)$.

We will check the following facts for the nonlinear mapping $\mathcal N$.

\textbf{ The mapping $\mathcal N$ is well-defined}. To this end, we need to show that for any $(F,G)\in \mathcal S_{\lambda,\mu}$, $\mathcal N(F,G)\in \mathcal S_{\lambda,\mu}$. We have from \eqref{eq:UniformLipschitzf_forw} and \eqref{eq:f_g_zero},
\begin{align*}
&\norme{\mathcal N(F,G)}_{\mathcal S_{\lambda,\mu}}^2=\esp\left(\int_{Q_T}\theta^{-2}\lambda^{-3}\mu^{-4}\xi^{-3}|f(y)|^2\dx\dt\right)+\esp\left(\int_{Q_T}\theta^{-2}\lambda^{-2}\mu^{-2}\xi^{-3}|g(y)|^2\dx\dt\right) \\
&\quad \leq \lambda^{-2}\mu^{-2} L^2 \esp\left(\int_{Q_T}\theta^{-2}\xi^{-3}|y|^2\dx\dt\right).
\end{align*}
At this point, we have also used that $\mu$ and $\lambda$ are much greater than one. 

Using \eqref{eq:est_weighted_spaces_forward} and $\norme{\xi^{-1}}_\infty\leq 1$ for all $(t,x)\in Q_T$, we get
\begin{align*}
&\norme{\mathcal N(F,G)}_{\mathcal S_{\lambda,\mu}}^2 \leq L^2 \lambda^{-2}\mu^{-2}\left(C_1\esp\norme{y_0}^2_{L^2(\dom)}+C\norme{(F,G)}_{\mathcal S_{\lambda,\mu}}^2 \right) <+\infty.
\end{align*}
This proves that $\mathcal N$ is well-defined.

\textbf{The mapping $\mathcal N$ is strictly contractive.} Let us consider couples $(F_i,G_i)\in  \mathcal S_{\lambda,\mu}$ for $i=1,2$. We denote the solutions of the corresponding equations by $y_1$ and $y_2$, respectively.  Using the fact that the nonlinearities $f$ and $g$ are globally Lipschitz, i.e. \eqref{eq:UniformLipschitzf_forw}, we have
\begin{align*}
\norme{\mathcal N(F_1,G_1)-\mathcal N(F_2,G_2)}^2_{\mathcal S_{\lambda,\mu}} &= \esp\left(\int_{Q_T}\theta^{-2}\lambda^{-3}\mu^{-4}\xi^{-3}|f(y_1)-f(y_2)|^2\dx\dt\right) \\
&\quad + \esp\left(\int_{Q_T}\theta^{-2}\lambda^{-2}\mu^{-2}\xi^{-3}|g(y_1)-g(y_2)|^2\dx\dt\right) \\
&\leq L^2\lambda^{-2}\mu^{-2}\esp\left(\int_{Q_T}\theta^{-2}|y_1-y_2|^2\dx\dt\right).
\end{align*}
Then applying \Cref{teo:contr_forward_source} and \Cref{rmk:uniquetrajectory}, and using the estimate \eqref{eq:est_weighted_spaces_forward} to the equation associated to $(F,G)=(F_1-F_2,G_1-G_2)$,  $y_0=0$, we deduce from the above inequality that
\begin{align}\notag 
\norme{\mathcal N(F_1,G_1)-\mathcal N(F_2,G_2)}^2_{\mathcal S_{\lambda,\mu}} &\leq C L^2\lambda^{-2}\mu^{-2} \left[\esp\left(\int_{Q_T}\theta^{-2}\lambda^{-3}\mu^{-4}\xi^{-3}|F_1-F_2|^2\dx\dt\right) \right. \\ \notag
&\hspace{2.6cm} \left.+\esp\left(\int_{Q_T}\theta^{-2}\lambda^{-2}\mu^{-2}\xi^{-3}|G_1-G_2|^2\dx\dt\right)\right]
\\ \label{eq:ineq_map}
&= C L^2 \lambda^{-2}\mu^{-2}\norme{(F_1,G_1)-(F_2,G_2)}^2_{\mathcal S_{\lambda,\mu}},
\end{align}
where $C=C(\dom,\dom_0)>0$ comes from \Cref{teo:contr_forward_source}. Observe that all the constants in the right-hand side of \eqref{eq:ineq_map} are uniform with respect to $\lambda$ and $\mu$ thus, if necessary, we can increase their value so $C L^2 \lambda^{-2}\mu^{-2}<1$. This yields that the mapping $\mathcal N$ is strictly contractive. 

Once we have verified these two conditions, by the Banach fixed point theorem, it follows that $\mathcal N$ has a unique fixed point $(F,G)$ in $\mathcal S_{\lambda,\mu}$. By setting $y$ the trajectory associated to the pair $(F,G)$, we observe that $y$ is the solution to \eqref{eq:forward_semilinear} and verifies $y(T,\cdot)=0$ in $\dom$, a.s. This concludes the proof of \Cref{th:semilinear_forward}.
\end{proof}

\section{Controllability of a semilinear backward stochastic parabolic equation}\label{sec:backward}


%

As for the forward equation, the main ingredient to prove \Cref{th:semilinear_backward} is a controllability result for a linear system with a source term. In this case, we shall focus on studying the controllability of 
\begin{equation*}
\begin{cases}
\d y=(-\Delta y+\chi_{\dom_0}h+F)\dt+Y\d{W}(t) &\text{in }Q_T, \\
y=0 &\text{on }\Sigma_T, \\
y(T)=y_T &\text{in }\dom,
\end{cases}
\end{equation*}
where $F\in L^2_{\mathcal F}(0,T;L^2(\dom))$ and $y_T\in L^2(\Omega,\mathcal F_T;L^2(\dom))$ are given. Unlike the previous section, we shall not devote to prove a Carleman estimate for the corresponding adjoint system (i.e. a forward equation). Although this is possible, we will see later that we can greatly simplify the problem by studying a random parabolic equation, for which a deterministic Carleman estimate will suffice. 

\subsection{A deterministic Carleman estimate and its consequence}

As we mentioned in \Cref{sec:new_carleman}, in \cite{BEG16} the authors have proved a Carleman estimate for the (backward) heat equation with weights that do not vanish as $t\to 0^{+}$ (see \eqref{eq:def_theta} and \eqref{eq:weights_0}). Following their approach it is possible to prove the analogous result for a forward equation. For this, we need to introduce some new weight functions which are actually the mirrored version of \eqref{eq:def_theta} and \eqref{eq:weights_0}. 

In more detail, let us consider the function $\beta$ as in \eqref{eq:prop_weight_beta} and let $0<T<1$. We define the function $\widetilde{\gamma}(t)$ as 
\begin{equation}\label{eq:def_theta_tilde}
\widetilde{\gamma}(t):=
\begin{cases}
\widetilde{\gamma}(t)=\frac{1}{t^m}, \quad t\in(0,T/4], \\
\widetilde{\gamma} \textnormal{ is decreasing on $[T/4,T/2]$}, \\
\widetilde{\gamma}(t)=1, \quad t\in[T/2,3T/4], \\
\widetilde{\gamma}(t)=1+\left(1-\frac{4(T-t)}{T}\right)^\sigma, \quad t\in[3T/4,T], \\
\widetilde{\gamma}\in C^2([0,T]),
\end{cases}
\end{equation}
where $m\geq 1$ and $\sigma\geq 2$ is defined on \eqref{def:sigma}. Observe that $\widetilde{\gamma}(t)$ is the mirrored version of $\gamma(t)$ in \eqref{eq:def_theta} with respect to $T/2$. Analogous to the properties of $\gamma$, the function $\widetilde\gamma$ preserves one important property which is that for the interval $[3T/4,T]$ the derivative of $\widetilde{\gamma}$ has a prescribed sign, i.e., $\widetilde{\gamma}_t\geq 0$. 

With this new function, we define the weights $\widetilde{\varphi}=\widetilde\varphi(t,x)$ and $\widetilde{\xi}=\widetilde{\xi}(t,x)$ as
\begin{equation}
\label{eq:weights_tilde}
\widetilde{\varphi}(t,x):=\widetilde{\gamma}(t)\left(e^{\mu(\beta(x)+6m)}-\mu e^{6\mu(m+1)}\right), \quad \widetilde{\xi}(t,x):=\tilde{\gamma}(t)e^{\mu(\beta(x)+6m)},
\end{equation}
where $\mu\geq 1$ is some parameter. In the same spirit, we set the weight $\widetilde{\theta}=\widetilde{\theta}(t,x)$ as 
\begin{equation*}
\widetilde{\theta}:=e^{\widetilde{\ell}} \quad\textnormal{where }\widetilde{\ell}(t,x)=\lambda \widetilde{\varphi}(t,x)
\end{equation*}
for a parameter $\lambda\geq 1$. 

In what follows, to keep the notation as light as possible and emphasizing that there is no possibility for confusion since the notation is specific for this section, we simply write $\widetilde\gamma=\gamma$, $\widetilde\theta=\theta$, and so on. 

We have the following Carleman estimate for the heat equation with source term
\begin{equation}
\label{eq:qsolheatSource}
\begin{cases}
\partial_t q-\Delta q=g(t,x) & \text{in } Q_T, \\
q= 0 &\text{on } \Sigma_T ,\\
q(0)=q_0(x) & \text{in } \dom.
\end{cases}
\end{equation}

\begin{theo}\label{car_refined}
For all $m\geq 1$, there exist constants $C>0$, $\lambda_0\geq 1$ and $\mu_0\geq 1$ such that for any $q_0\in L^2(\dom)$ and any $g\in L^2(Q_T)$, the weak solution to \eqref{eq:qsolheatSource}
satisfies
\begin{align*}
&\int_{Q_T}\theta^2\lambda\mu^2\xi|\nabla q|^2\dx\dt+\int_{Q_T}\theta^2\lambda^3\mu^4\xi^{3}|q|^2\dx\dt\notag\\
& + \int_{\dom} \theta^2(T)|\nabla q(T)|^2\dx +  \int_{\dom} \lambda^2\mu^3e^{2\mu(6m+1)}\theta^2(T)|q(T)|^2 \dx\notag \\
& \leq C\left(\int_{Q_T} \theta^2|g|^2\dx\dt+\iint_{\mathcal \dom_0\times(0,T)}\theta^2\lambda^3\mu^4\xi^3 |q|^2\dx\dt \right),
\end{align*}%
for all $\mu\geq \mu_0$ and $\lambda\geq \lambda_0$.
\end{theo}

The proof of this result is a straightforward adaptation of \cite[Theorem 2.5]{BEG16}, just by taking into account that in this case the weight $\gamma$ verifies $\gamma_t\geq 0$ in $[3T/4,T]$, contrasting with the fact that $\gamma_t\leq 0$ in [0,T/4] as in \cite{BEG16} or as we have used in the proof of \Cref{thm:carleman_backward_0}.

Let us consider the forward parabolic equation given by

\begin{equation}\label{eq:forw_gen}
\begin{cases}
\d q=(\Delta q+G_1)\dt+ G_2\d{W}(t) & \textnormal{in } Q_T, \\
q= 0 &\textnormal{on } \Sigma_T ,\\
q(0,x)=q_0(x) & \textnormal{in } \dom,
\end{cases}
\end{equation}
where $G_i\in L^2_{\mathcal F}(0,T;L^2(\dom))$, $i=1,2$, and $q_0\in L^2(\Omega,\mathcal F_0;L^2(\dom))$. An immediate consequence of \Cref{car_refined} is a Carleman estimate for a random parabolic equation. More precisely, we have the following.

\begin{lem} \label{lem:car_random}
Assume that $G_2\equiv 0$. For all $m\geq 1$, there exists constants $C>0$, $\lambda_0\geq 1$ and $\mu_0\geq 1$ such that for any $q_0\in L^2(\Omega,\mathcal F_0; L^2(\dom))$ and any $g\in L^2_{\mathcal F}(0,T;L^2(\dom))$, the corresponding solution to \eqref{eq:forw_gen} with $G_2=0$ satisfies
\begin{align}
&\esp\left(\int_{Q_T}\theta^2\lambda\mu^2\xi|\nabla q|^2\dx\dt\right)+\esp\left(\int_{Q_T}\theta^2\lambda^3\mu^4\xi^{3}|q|^2\dx\dt\right)\notag\\
& + \esp\left(\int_{\dom} \theta^2(T)|\nabla q(T)|^2\dx\right) +  \esp\left(\int_{\dom} \lambda^2\mu^3e^{2\mu(6m+1)}\theta^2(T)|q(T)|^2 \dx\right) \notag \\ \label{car_sigma}
& \leq C\esp\left(\int_{Q_T} \theta^2|g|^2\dx\dt+\iint_{\mathcal \dom_0\times(0,T)}\theta^2\lambda^3\mu^4\xi^3 |q|^2\dx\dt \right),
\end{align}%
for all $\mu\geq \mu_0$ and $\lambda\geq \lambda_0$.
\end{lem}

\subsection{A controllability result for a linear backward stochastic heat equation with source term and one control}

Inspired by the duality technique presented in \cite[Prop. 2.2]{liu14}, we present a controllability result for a linear backward stochastic heat equation with a source term. To this end, consider the linear control system given by %
\begin{equation}\label{eq:backward_source}
\begin{cases}
\d y=(-\Delta y+\chi_{\dom_0}h+F)\dt+Y\d{W}(t) &\text{in }Q_T, \\
y=0 &\text{on }\Sigma_T, \\
y(T)=y_T &\text{in }\dom, 
\end{cases}
\end{equation}
where $F\in L^2_{\mathcal F}(0,T;L^2(\dom))$ is a given fixed source term and $h\in L^2_{\mathcal F}(0,T;L^2(\dom_0))$ is a control. 

 In what follows, we consider constants $\mu$ and $\lambda$ large enough such that \eqref{car_sigma} holds. We define the space $\widetilde{\mathcal{S}}_{\lambda,\mu}:=\left\{F\in L^2_{\mathcal F}(0,T;L^2(\dom)):\esp\left(\int_{Q_T}\theta^{-2}\lambda^{-3}\mu^{-4}\xi^{-3}|F|^2\dx\dt\right)<+\infty\right\}$, endowed with the canonical norm. We have the following global null-controllability result for system \eqref{eq:backward_source}.
\begin{theo}\label{teo:contr_backward_source}
For any initial datum $y_T\in L^2(\Omega,\mathcal F_T;L^2(\dom))$ and any $F\in \widetilde{\mathcal{S}}_{\lambda,\mu}$, there exists a control $\widehat{h}\in L^2(0,T;L^2(\dom_0))$ such that the associated solution $(\widehat{y},\widehat{Y})\in [L^2_{\mathcal F}(\Omega;C[0,T];L^2(\dom))\cap L^2_{\mathcal F}(0,T;H_0^1(\dom))]\times L^2_{\mathcal F}(0,T;L^2(\dom))$ to system \eqref{eq:backward_source} satisfies $\widehat{y}(0)=0$ in $\dom$, a.s. Moreover, the following estimate holds
\begin{align}\notag 
\esp&\left(\int_{Q_T}\theta^{-2}|\widehat{y}|^2\dx\dt\right)+\esp\left(\int_{Q_T}\theta^{-2}\lambda^{-2}{\mu^{-2}}\xi^{-2}|\widehat{Y}|^2\dx\dt\right)+\esp\left(\int_{Q_T}\theta^{-2}\lambda^{-3}\mu^{-4}\xi^{-3}|\widehat{h}|^2\dx\dt\right) \\ \label{eq:est_control_weighted_spaces}
&\quad \leq C_1\esp\left(\|y_T\|^2_{L^2(\dom)}\right)+C\esp\left(\int_{Q_T}\theta^{-2}\lambda^{-3}\mu^{-4}\xi^{-3}|F|^2\dx\dt\right),
\end{align}
where $C_1>0$ is a constant depending on $\dom,\dom_0,\mu,\lambda$ and $C>0$ only depends on $\dom$ and $\dom_0$. 
\end{theo}
\begin{rmk}
\label{rmk:uniquetrajectory_backward}
As before, from classical arguments, see e.g. \cite[Proposition 2.9]{LLT13}, from \Cref{teo:contr_backward_source} we can construct a linear continuous mapping that associates every initial datum $y_T\in L^2(\Omega,\mathcal F_T;L^2(\dom))$ and every source term $F\in \widetilde{\mathcal{S}}_{\lambda,\mu}$, to a trajectory $(\widehat{y}, \widehat{h})$ such that $\widehat{y}(0)=0$ in $\dom$, a.s. and \eqref{eq:est_control_weighted_spaces} holds.
\end{rmk}
\begin{proof} The proof is very similar to the one of \Cref{teo:contr_forward_source} and requires only some adaptations. We emphasize their main differences. 

\smallskip
\textbf{Step 1}. For any $\epsilon>0$, let us consider the weight function $\gamma_\epsilon(t)$ given by
\begin{equation*}
\gamma_\epsilon(t):=
\begin{cases}
\gamma_\epsilon(t)=\gamma(t+\epsilon), \quad t\in [0,T/2-\epsilon], \\
\gamma_\epsilon(t)=1, \quad t\in[T/2-\epsilon,3T/4] \\
\gamma_\epsilon(t)=1+(1+\frac{4(T-t)}{T})^{\sigma}, \quad t\in [3T/4,T], \\
\sigma \textnormal{ as in \eqref{def:sigma}}. 
\end{cases}
\end{equation*}
In this way, $\gamma_\epsilon(t)\leq \gamma(t)$ for $t\in[0,T]$. We set the corresponding weight  $\varphi_\epsilon$ as in \eqref{eq:weights_tilde} by replacing the function $\gamma$ by $\gamma_\epsilon$. Also, we write $\theta_\epsilon=e^{\lambda\varphi_\epsilon}$. 

We introduce the cost functional
\begin{align*}\notag
\mathcal I_\epsilon(h):=&\frac{1}{2}\esp\left(\int_{Q_T}\theta_\epsilon^{-2}|y|^2\dx\dt\right)+\frac{1}{2}\esp\left(\int_{0}^T\!\!\!\int_{\dom_0}\theta^{-2}\lambda^{-3}\mu^{-4}\xi^{-3}|h|^2\dx\dt\right)\\ \label{eq:func}
&+\frac{1}{2\epsilon}\esp\left(\int_{\dom}|y(0)|^2\dx\right)
\end{align*}
and consider the minimization problem
\begin{equation}\label{eq:prob_min_forw}
\begin{cases}
\min_{h\in \mathcal H} \mathcal I_\epsilon(h), \\
\textnormal{subject to equation \eqref{eq:backward_source},} 
\end{cases}
\end{equation}
where \begin{equation*}
\mathcal H=\left\{h\in L^2_{\mathcal F}(0,T;L^2(\dom_0)):\esp\left(\int_{0}^{T}\!\!\!\int_{\dom_0}\theta^{-2}\lambda^{-3}\mu^{-4}\xi^{-3}|h|^2\dx\dt\right)<+\infty\right\}.
\end{equation*}

It can be readily seen that the functional $\mathcal I_\epsilon$ is continuous, strictly convex and coercive. Therefore, the minimization problem \eqref{eq:prob_min_forw} admits a unique optimal solution that we denote by $h_\epsilon$. As in the proof of \Cref{teo:contr_forward_source} the minimizer $h_\epsilon$ can be characterized as 
\begin{equation}\label{eq:h_epsilon_forw}
h_\epsilon=\chi_{\dom_0}\lambda^3\mu^{4}\xi^3\theta^2q_\epsilon \quad\text{in Q}, \quad a.s.,
\end{equation}
where $q_\epsilon$ verifies the random forward equation
\begin{equation}\label{eq:q_optim_forw}
\begin{cases}
\d{q_\epsilon}=(\Delta q_\epsilon+\theta^{-2}_\epsilon y_\epsilon)\dt &\text{in }Q_T, \\
q_\epsilon=0 &\text{on }\Sigma_T, \\
q_\epsilon(0)=\frac{1}{\epsilon}y_\epsilon(0) &\text{in }\dom,
\end{cases}
\end{equation}
and where $(y_\epsilon,y_\epsilon(0))$ can be extracted from $(y_\epsilon,Y_\epsilon)$ the solution to \eqref{eq:backward_source} with control $h=h_\epsilon$. Observe that since $y_\epsilon\in L^2_{\mathcal F}(\Omega; C([0,T];L^2(\dom)))$ the evaluation of $y_\epsilon$ at $t=0$ is meaningful and \eqref{eq:q_optim_forw} is well-posed for any $\epsilon>0$. Also notice that there is no term containing $W(t)$ so \eqref{eq:q_optim_forw} is regarded as a random equation. This greatly simplifies our task, since we only need to use the Carleman estimate of \Cref{lem:car_random} to deduce the uniform estimate for the solutions to $(y_\epsilon,Y_\epsilon)$ in the next step.  

\smallskip
\textbf{Step 2.} Using It\^{o}'s formula, we can compute $\d(y_\epsilon q_\epsilon)$ and deduce
\begin{align*}
\esp\left(\int_{\dom}y_\epsilon(T)q_\epsilon(T)\dx\right)&=\esp\left(\int_{\dom}y_\epsilon(0)q_\epsilon(0)\dx\right)+\esp\left(\int_{Q}(-\Delta y_\epsilon+F+\chi_{\dom_0}h_\epsilon)q_\epsilon\dx\dt\right) \\
&\quad + \esp\left(\int_{Q}(\Delta q_\epsilon+\theta^{-2}_\epsilon y_\epsilon)y_\epsilon\dx\dt\right)
\end{align*}
whence, using equations \eqref{eq:backward_source}, \eqref{eq:q_optim_forw}, and identity \eqref{eq:h_epsilon_forw}, we get
\begin{align}\notag
\esp&\left(\int_{0}^{T}\!\!\!\int_{\dom_0}\lambda^3\mu^4\xi^3\theta^2|q_\epsilon|^2\dx\dt\right)+\esp\left(\int_{Q_T}\theta^{-2}_\epsilon|y_\epsilon|^2\dx\dt\right)+\frac{1}{\epsilon}\esp\left(\int_{\dom}|y_\epsilon(0)|^2\dx\right) \\ \label{eq:id_prod_yq_forw}
&=\esp\left(\int_{\dom}y_Tq_\epsilon(T)\dx\right)-\esp\left(\int_{Q_T}F q_\epsilon \dx\dt \right).
\end{align}

In view of \eqref{car_sigma}, we use Cauchy-Schwarz and Young inequalities in the right-hand side of \eqref{eq:id_prod_yq_forw} to introduce the weight function as follows
\begin{align}\notag
\esp&\left(\int_{0}^{T}\!\!\!\int_{\dom_0}\theta^2 \lambda^3\mu^4\xi^3e^{-2s\varphi}|q_\epsilon|^2\dx\dt\right)+\esp\left(\int_{Q_T}\theta_\epsilon^{-2}|y_\epsilon|^2\dx\dt\right)+\frac{1}{\epsilon}\esp\left(\int_{\dom}|y_\epsilon(0)|^2\dx\right) \\ \notag
&\leq \delta\left[\esp\left(\int_{\dom}\lambda^2\mu^3\theta^{2}(T)|q_\epsilon(T,x)|^2\dx\right)+\esp\left(\int_{Q_T}\theta^2 \lambda^3\mu^4\xi^3|q_\epsilon|^2\dx\dt\right)\right]\\ \label{eq:des_prod_yq_forw}
&\quad + C_\delta\left[\esp\left(\int_{\dom}\lambda^{-2}\mu^{-3}\theta^{-2}(T)|y_T|^2\dx\right)+\esp\left(\int_{Q_{T}}\theta^{-2}\lambda^{-3}\mu^{-4}\xi^{-3}|F|^2\dx\dt\right)\right]
\end{align}
with $\delta>0$. Applying inequality \eqref{car_sigma} to \eqref{eq:q_optim_forw} and using it to estimate in the right-hand side of \eqref{eq:des_prod_yq_forw}, we obtain, after taking $\delta>0$ small enough, that
\begin{align}\notag
\esp&\left(\int_{0}^{T}\!\!\!\int_{\dom_0}\theta^{2}\lambda^3\mu^4\xi^3|q_\epsilon|^2\dx\dt\right)+\esp\left(\int_{Q_T}\theta_{\epsilon}^{-2}|y_\epsilon|^2\dx\dt\right)+\frac{1}{\epsilon}\esp\left(\int_{\dom}|y_\epsilon(0)|^2\dx\right) \\ \notag
&\leq C\left[\esp\left(\int_{\dom}\lambda^{-2}\mu^{-3}\theta^{-2}(T)|y_T|^2\dx\right)+\esp\left(\int_{Q_{T}}\theta^{-2}\lambda^{-3}\mu^{-4}\xi^{-3}|F|^2\dx\dt\right)\right]
\end{align}
for some constant $C>0$ only depending on $\dom$ and $\dom_0$. At this point, we have used the fact that $\theta^2 \theta_\epsilon^{-2}\leq 1$ for all $(t,x)\in Q_T$.

Recalling the characterization of the optimal control $h_\epsilon$ in \eqref{eq:h_epsilon_forw} we obtain
\begin{align}\notag
\esp&\left(\int_{0}^{T}\!\!\!\int_{\dom_0}\theta^{-2}\lambda^{-3}\mu^{-4}\xi^{-3}|h_\epsilon|^2\dx\dt\right)+\esp\left(\int_{Q_T}e^{2s\varphi_\epsilon}|y_\epsilon|^2\dx\dt\right)+\frac{1}{\epsilon}\esp\left(\int_{\dom}|y_\epsilon(0)|^2\dx\right) \\ \label{eq:iden_uniform_forw}
&\leq C\left[\esp\left(\int_{\dom}\theta^{-2}\lambda^{-2}\mu^{-3}|y_T|^2\dx\right)+\esp\left(\int_{Q_{T}}\theta^{-2}\lambda^{-3}\mu^{-4}\xi^{-3}|F|^2s^{-3}\xi^{-3}\dx\dt\right)\right].
\end{align}

Now, our task is to add a weighted integral of the process $Y$ on the left-hand side of the above inequality. To do that, using It\^{o}'s formula and equation \eqref{eq:backward_source} with $h=h_\epsilon$ yield
\begin{align*}
\d(\theta^{-2}_\epsilon \lambda^{-2}\xi^{-2}y_\epsilon^2)=&(\theta^{-2}_\epsilon \lambda^{-2}\xi^{-2})_ty_\epsilon^2\dt+\theta^{-2}_\epsilon \lambda^{-2}\xi^{-2}Y_\epsilon^2\\
&+2\theta^{-2}_\epsilon \lambda^{-2}\xi^{-2} y_\epsilon\left[(-\Delta y_\epsilon+\chi_{\dom_0}h_\epsilon+F)\dt+Y_\epsilon\d{W}(t)\right]
\end{align*}
 and after some integration by parts and substituting the initial datum, we get
\begin{align}\notag 
\esp&\left(\int_{Q_T}\theta^{-2}_\epsilon \lambda^{-2}\xi^{-2}|Y_\epsilon|^2\dx\dt\right)+2\esp\left(\int_{Q_T}\theta_\epsilon^{-2}\lambda^{-2}\xi^{-2}|\nabla y_\epsilon|^2\dx\dt\right)\\ \notag
&+\esp\left(\int_{Q_T}(\theta_\epsilon ^{-2}\lambda^{-2}\xi^{-2})_t|y_\epsilon|^2\dx\dt\right)=\esp\left(\int_{\dom}\theta_\epsilon^{-2}(T)\lambda^{-2}\xi^{-2}(T)|y_T|^2\dx\right) \\ \notag
&-2\esp\left(\int_{Q_T}\nabla(\theta_\epsilon^{-2} \lambda^{-2}\xi^{-2})\cdot\nabla y_{\epsilon}y_\epsilon\dx\dt\right)-2\esp\left(\int_{Q_T}\theta_\epsilon^{-2}\lambda^{-2}\xi^{-2}y_\epsilon F\dx\dt\right) \\ \label{eq:iden_Y_weight_forw}
&-2\esp\left(\int_{0}^{T}\!\!\!\int_{\dom_0}\theta_\epsilon^{-2}\lambda^{-2}\xi^{-2}y_\epsilon h_\epsilon \dx\dt\right)
\end{align}
Observe that the term containing $y_T$ is well defined since, by construction, the weight $\theta_\epsilon^{-1}$ does not blow up at $t=T$. Also, notice that there is no term of $y_\epsilon(0,x)$ since $\xi^{-1}(0)=0$ and the weight $\theta_\epsilon^{-1}$ does not blow up at $t=0$.

Let us analyze the term containing $(\theta_\epsilon^{-2}\lambda^{-2}\xi^{-2})_t$ in the left-hand side of the above identity. We split the integral as
\begin{align}\notag
\esp\left(\int_{Q_T}(\theta^{-2}_\epsilon \lambda^{-2}\xi^{-2})_t|y_\epsilon|^2\dx\dt\right)&=\esp\left(\int_{0}^{3T/4}\!\!\!\int_{\dom}(\theta_\epsilon^{-2}\lambda^{-2}\xi^{-2})_t|y_\epsilon|^2\dx\dt\right)\\ \label{eq:split_deriv}
&\quad +\esp\left(\int_{3T/4}^{T}\!\int_{\dom}(\theta_\epsilon^{-2}\lambda^{-2}\xi^{-2})_t|y_\epsilon|^2\dx\dt\right).
\end{align}
We note that for $t\in[3T/4,T]$, $\gamma_\epsilon(t)=\gamma(t)$, so we can drop the dependence of $\epsilon$. Also notice that on this time interval $\gamma_t\geq 0$ and $1\leq \gamma \leq 2$. Thus, computing explicitly, we have
\begin{equation}\label{eq:iden_deriv_eps}
(\theta^{-2}\lambda^{-2}\xi^{-2})_t=-2\theta^{-2}\lambda^{-1}\frac{\gamma_t}{\gamma}\varphi \xi^{-2}-2\theta^{-2}\lambda^{-2}\frac{\gamma_t}{\gamma}\xi^{-2}.
\end{equation}
Recall that $\varphi<0$, thus
\begin{equation}
(\theta^{-2}\lambda^{-2}\xi^{-2})_t \geq c \theta^{-2} \lambda^{-1} \gamma_t |\varphi| \xi^{-2}.
\end{equation}
for all $t\in[3T/4,T]$, where $c>0$ only depends on $\dom$ and $\dom_0$. Therefore, 
\begin{equation}\label{eq:est_deriv_34_T}
\esp\left(\int_{3T/4}^{T}\!\int_{\dom}(\theta_\epsilon^{-2}\lambda^{-2}\xi^{-2})_t|y_\epsilon|^2\dx\dt\right)\geq 0
\end{equation}
and this term can be dropped. For $t\in[0,3T/4]$, we can use expression \eqref{eq:iden_deriv_eps} (replacing everywhere the weights depending on $\epsilon$) and the fact that $|\partial_t\gamma_\epsilon|\leq C\gamma_\epsilon^2$ to obtain
\begin{equation*}
\abs{(\theta_\epsilon^{-2}\lambda^{-2}\xi^{-2})_t}\leq C\theta_\epsilon^{-2}\lambda^{-1}\mu,
\end{equation*}
where the constant $C>0$ is uniform with respect to $\lambda$ and $\mu$. Therefore,
\begin{equation}\label{eq:est_0_34}
\abs{\esp\left(\int_{0}^{3T/4}\!\!\!\int_{\dom}(\theta_\epsilon^{-2}\lambda^{-2}\xi^{-2})_t|y_\epsilon|^2\dx\dt\right)} \leq  C\esp\left(\int_{0}^{3T/4}\!\!\!\int_{\dom}\theta_\epsilon^{-2}\lambda^{-1}\mu|y_\epsilon|^2\dx\dt\right).
\end{equation}
Thus, using formulas \eqref{eq:split_deriv} and \eqref{eq:est_deriv_34_T}--\eqref{eq:est_0_34} we deduce from \eqref{eq:iden_Y_weight_forw} that
\begin{align}\notag 
\esp&\left(\int_{Q_T}\theta_\epsilon^{-2}\lambda^{-2}\xi^{-2}|Y_\epsilon|^2\dx\dt\right)+2\esp\left(\int_{Q_T}\theta_\epsilon^{-2}\lambda^{-2}\xi^{-2}|\nabla y_\epsilon|^2\dx\dt\right)\\ \notag
&\leq C\esp\left(\int_{\dom}\theta^{-2}(T) \lambda^{-2}|y_T|^2\dx\right)+C\esp\left(\int_{0}^{3T/4}\!\!\!\int_{\dom}\theta_\epsilon^{-2}\lambda^{-1}\mu |y_\epsilon|^2\dx\dt\right)\\ \notag
&+2\abs{\esp\left(\int_{Q_T}\nabla(\theta_\epsilon^{-2}\lambda^{-2}\xi^{-2})\cdot\nabla y_{\epsilon}y_\epsilon\dx\dt\right)}+2\abs{\esp\left(\int_{Q_T}\theta_\epsilon^{-2}\lambda^{-2}\xi^{-2}y_\epsilon F\dx\dt\right)} \\ \label{eq:ineq_Y_weight_forw}
&+2\abs{\esp\left(\int_{0}^{T}\!\!\!\int_{\dom_0}\theta_\epsilon^{-2}\lambda^{-2}\xi^{-2}y_\epsilon h_\epsilon \dx\dt\right)}.
\end{align}
For the first term in the right-hand side, we have used that $\theta^{-2}_\epsilon(T)=\theta^{-2}(T)$ and $\xi^{-1}(T)\leq C$ for some $C>0$ only depending on $\dom$ and $\dom_0$.

Employing Cauchy-Schwarz and Young inequalities, we estimate the last three terms of the above inequality. For the first one, we have
\begin{align}\notag
2&\abs{\esp\left(\int_{Q_T}\nabla(\theta_\epsilon^{-2}\lambda^{-2}\xi^{-2})\cdot\nabla y_{\epsilon}y_\epsilon\dx\dt\right)}\\ \label{eq:first_est_forw}
&\leq \delta \esp\left(\int_{Q_T}\theta_\epsilon^{-2}\lambda^{-2}\xi^{-2}|\nabla y_\epsilon|^2\dx\dt\right) + C(\delta)\esp\left(\int_{Q_T}\theta_\epsilon^{-2}\mu^2|y_\epsilon|^2\dx\dt\right)
\end{align}
for any $\delta>0$. Here, we have used that $|\nabla(\theta_\epsilon^{-2}\lambda^{-2}\xi^{-2})|\leq C\theta_\epsilon^{-2}\mu \lambda^{-1} \xi^{-1}$. For the second term, we get
\begin{align}\notag 
2&\abs{\esp\left(\int_{Q_T}\theta_\epsilon^{-2}\lambda^{-2}\xi^{-2}y_\epsilon F\dx\dt\right)} \\ \label{eq:second_est_forw}
&\leq \esp\left(\int_{Q_T}\theta_\epsilon^{-2}\mu^2\lambda^{-1}\xi^{-1}|y_\epsilon|^2\dx\dt\right)+\esp\left(\int_{Q_T}\theta^{-2}\lambda^{-3}\mu^{-2}\xi^{-3}|F|^2\dx\dt\right),
\end{align}
where we have used that $\theta_\epsilon^{-2}\leq \theta^{-2}$. For the last one, we readily have
\begin{align}\notag
2&\abs{\esp\left(\int_{0}^{T}\!\!\!\int_{\dom_0}\theta_\epsilon^{-2}\lambda^{-2}\xi^{-2}y_\epsilon h_\epsilon \dx\dt\right)}\\\label{eq:third_est_forw}
&\leq \esp\left(\int_{Q_T}\theta_\epsilon^{-2}\mu^2|y_\epsilon|^2\dx\dt\right)+\esp\left(\int_{0}^{T}\!\!\!\int_{\dom_0}\theta_\epsilon^{-2}\lambda^{-4}\mu^{-2}\xi^{-4}|h_\epsilon|^2\dx\dt\right).
\end{align}

Using estimates \eqref{eq:first_est_forw}--\eqref{eq:third_est_forw} in \eqref{eq:ineq_Y_weight_forw} and taking $\delta>0$ small enough, we deduce after collecting similar terms that
\begin{align}\notag 
\esp&\left(\int_{Q_T}\theta_\epsilon^{-2}\lambda^{-2}\mu^{-2}\xi^{-2}|Y_\epsilon|^2\dx\dt\right)+\esp\left(\int_{Q_T}\theta_\epsilon^{-2}\lambda^{-2}\mu^{-2}\xi^{-2}|\nabla y_\epsilon|^2\dx\dt\right) \\ \notag
&\leq C\esp\left(\int_{\dom}\theta^{-2}(T)\lambda^{-2}\mu^{-2}|y_T|^2\dx\right)+C\esp\left(\int_{Q_T}\theta_\epsilon^{-2}|y_\epsilon|^2\dx\dt\right)\\\label{eq:ineq_Y_weight_final_forw}
&\quad +C\esp\left(\int_{Q_T}\theta^{-2}\lambda^{-3}\mu^{-4}\xi^{-3}|F|^2\dx\dt\right) +C\esp\left(\int_{0}^{T}\!\!\!\int_{\dom_0}\theta_\epsilon^{-2}\lambda^{-3}\mu^{-4}\xi^{-3}|h_\epsilon|^2 \dx\dt\right).
\end{align}
At this point, we have adjusted the powers of $\lambda$ and $\xi$ in the last term by using the fact that $\lambda^{-1}\xi^{-1}\leq C$ for some constant only depending on $\dom,\dom_0$. 

Finally, combining \eqref{eq:ineq_Y_weight_final_forw} and \eqref{eq:iden_uniform_forw} we get
\begin{align}\notag
\esp&\left(\int_{0}^{T}\!\!\!\int_{\dom_0}\theta_\epsilon^{-2}\lambda^{-3}\mu^{-4}\xi^{-3}|h_\epsilon|^2\dx\dt\right)+\esp\left(\int_{Q_T}\theta_\epsilon^{-2}(|y_\epsilon|^2+\lambda^{-2}\mu^{-2}\xi^{-2}|Y_\epsilon|^2)\dx\dt\right) \\ \notag
&+\esp\left(\int_{Q_T}\theta_\epsilon^{-2}\lambda^{-2}\mu^{-2}\xi^{-2}|\nabla y_\epsilon|^2\dx\dt\right)+\frac{1}{\epsilon}\esp\left(\int_{\dom}|y_\epsilon(0)|^2\dx\right)  \\ \label{eq:iden_uniform_final_forw}
&\leq C\left[\esp\left(\int_{\dom}\theta^{-2}(T)\lambda^{-2}\mu^{-2}|y_T|^2\dx\right)+\esp\left(\int_{Q_{T}}\theta^{-2}\lambda^{-3}\mu^{-4}\xi^{-3}|F|^2 \dx\dt\right)\right],
\end{align}
for some positive constant $C$ only depending on $\dom,\dom_0$.

\textbf{Step 3.} The last step is essentially the same as in the proof of \Cref{teo:contr_forward_source}. Since the right-hand side of \eqref{eq:iden_uniform_final_forw} is uniform with respect to $\epsilon$, we readily deduce that there exists $(\widehat{h},\widehat{y},\widehat{Y})$ such that
\begin{equation}\label{eq:weak_conv_forw}
\begin{cases}
h_\epsilon\weakly \widehat{h} &\textnormal{weakly in } L^2(\Omega\times(0,T);L^2(\dom_0)), \\
y_{\epsilon}\weakly \widehat{y} &\textnormal{weakly in } L^2(\Omega\times(0,T);H_0^1(\dom)), \\
Y_\epsilon\weakly Y &\textnormal{weakly in } L^2(\Omega\times(0,T);L^2(\dom)).
\end{cases}
\end{equation}
To check that $(\widehat{y},\widehat{Y})$ is the solution to \eqref{eq:backward_source} associated to $\widehat{h}$ can be done exactly as in \Cref{teo:contr_forward_source}, so we omit it. 

To conclude, we notice from \eqref{eq:iden_uniform_forw} that $\widehat{y}(0)=0$ in $\dom$, a.s. Also, from the weak convergence \eqref{eq:weak_conv_forw}, Fatou's lemma and the uniform estimate \eqref{eq:iden_uniform_forw} we deduce \eqref{eq:est_control_weighted_spaces}. This ends the proof of \Cref{teo:contr_backward_source}.
\end{proof}

\subsection{Proof the nonlinear result for the backward equation}

Now, we are in position to prove \Cref{th:semilinear_backward}. The proof is very similar to the one of \Cref{th:semilinear_forward} but for the sake of completeness, we give it. 

Let us fix the parameters $\lambda$ and $\mu$ in \Cref{teo:contr_backward_source} to a fixed value sufficiently large. Recall that in turn, this parameter comes from the \Cref{car_refined} and should be selected as $\lambda\geq \lambda_0$ and $\mu\geq \mu_0$ for some $\lambda_0,\mu_0\geq 1$, so there is no contradiction. 

Let us consider a nonlinearity $f$ fulfilling \eqref{eq:UniformLipschitzf} and \eqref{eq:fzero} and define
\begin{equation*}
\widetilde{\mathcal N}: F\in \mathcal{\widetilde{S}}_{\lambda,\mu}\mapsto f(\omega,t,x,y,Y)\in \mathcal{\widetilde{S}}_{\lambda,\mu},
\end{equation*}
where $(y,Y)$ is the trajectory of \eqref{eq:backward_source} associated to the data $y_T$ and $F$, defined by \Cref{teo:contr_backward_source} and \Cref{rmk:uniquetrajectory_backward}.  In what follows, to abridge the notation, we simply write $f(y,Y)$.

We will check the following facts for the nonlinear mapping $\widetilde{\mathcal N}$.

\textbf{The mapping $\mathcal N$ is well-defined}. To this end, we need to show that for any $F\in \mathcal{\widetilde{S}}_{\lambda,\mu}$, $\mathcal N(F)\in \mathcal{\widetilde{S}}_{\lambda,\mu}$. We have from \eqref{eq:UniformLipschitzf} and \eqref{eq:fzero}
\begin{align*}
\norme{\widetilde{\mathcal N}(F)}_{\mathcal{\widetilde{S}}_{\lambda,\mu}}^2&=\esp\left(\int_{Q_T}\theta^{-2}\lambda^{-3}\mu^{-4}\xi^{-3}|f(y,Y)|^2\dx\dt\right) \\
&\leq 2 L^2\esp\left(\int_{Q_T}\theta^{-2}\lambda^{-3}\mu^{-4}\xi^{-3}\left[|y|^2+|Y|^2\right]\dx\dt\right) \\
&\leq 2 L^2 \lambda^{-1}\mu^{-2} \left[ \esp\left(\int_{Q_T}\theta^{-2}\lambda^{-2}\mu^{-2}\xi^{-2}|Y|^2\dx\dt\right)+\esp\left(\int_{Q_T}\theta^{-2}|y|^2\dx\dt\right)\right] \\
&\leq 2L^2 \lambda^{-1}\mu^{-2}\left(C_1\esp\left[\norme{y_T}^2_{L^2(\dom)}\right)+C\esp\left(\int_{Q_T}\theta^{-2}\lambda^{-3}\mu^{-4}\xi^{-3}|F|^2\dx\dt\right)\right]\\
&<+\infty,
\end{align*}
where we have used \eqref{eq:est_control_weighted_spaces} and that $\norme{\xi^{-1}}_\infty\leq 1$.  This proves that $\mathcal N$ is well-defined.

\textbf{The mapping $\mathcal N$ is a strictly contraction mapping.} Let us consider $F_i\in \mathcal \mathcal{\widetilde{S}}_{\lambda,\mu}$, $i=1,2$.  From the properties of the nonlinearity $f$, we have
\begin{align*}\notag 
&\norme{\widetilde{\mathcal N}(F_1)-\widetilde{\mathcal N}(F_2)}^2_{\mathcal{\widetilde{S}}_{\lambda,\mu}}\\
&\quad =\esp\left(\int_{Q_T}\theta^{-2}\lambda^{-3}\mu^{-4}\xi^{-3}|f(y_1,Y_1)-f(y_2,Y_2)|^2\dx\dt\right) \\ \notag
&\quad \leq 2L^2\lambda^{-1}\mu^{-2} \esp\left(\int_{Q_T}\theta^{-2}\lambda^{-2}\mu^{-2}\xi^{-2}|Y_1-Y_2|^2\dx\dt+\int_{Q_T}\theta^{-2}|y_1-y_2|^2\dx\dt\right).
\end{align*}
Then applying \Cref{teo:contr_backward_source} to the equation associated to $F=F_1-F_2$,  $y_T=0$, and using the corresponding estimate \eqref{eq:est_control_weighted_spaces}, we deduce from the above inequality that
\begin{align}\notag 
\norme{\widetilde{\mathcal N}(F_1)-\widetilde{\mathcal N}(F_2)}^2_{\mathcal{\widetilde{S}}_{\lambda,\mu}} &\leq 2C  L^2\lambda^{-1}\mu^{-2} \esp\left(\int_{Q_T}\theta^{-2}\lambda^{-3}\mu^{-4}\xi^{-3}|F_1-F_2|^2\dx\dt\right) \\ \label{eq:ineq_map_back}
&= 2CL^2 \lambda^{-1}\mu^{-2}\norme{F_1-F_2}^2_{\mathcal{\widetilde{S}}_{\lambda,\mu}},
\end{align}
where $C=C(\dom,\dom_0)>0$ comes from \Cref{teo:contr_backward_source}. Observe that all the constants in the right-hand side of \eqref{eq:ineq_map_back} are uniform with respect to $\lambda$ and $\mu$ thus, if necessary, we can increase the value of $\lambda$ and $\mu$ so $CL^2 \lambda^{-1}\mu^{-2}<1$. This yields that the mapping is strictly contractive. 

Once we have verified these two conditions, it follows that $\widetilde{\mathcal N}$ has a unique fixed point $F$ in $\mathcal{\widetilde{S}}_{\lambda,\mu}$. By setting $(y,Y)$ the trajectory associated to this $F$, we observe that $(y,Y)$ is the solution to \eqref{eq:backward_semilinear} and verifies $y(0,\cdot)=0$ in $\dom$, a.s. This concludes the proof of \Cref{th:semilinear_backward}.

\section{Further results and remarks}\label{sec:conclusion}

\subsection{A new Carleman estimate for forward equation as a consequence of \Cref{teo:contr_backward_source}}

The controllability result provided by \Cref{teo:contr_backward_source} yields as a byproduct the obtention of a new global Carleman estimate for forward stochastic parabolic equations with a weight that do not vanish as $t\to T^-$. In fact, under the construction of weights shown in \eqref{eq:def_theta_tilde} and \eqref{eq:weights_tilde} (where again we drop the tilde notation for simplicity), we are able to prove the following result.
\begin{prop}\label{prop:car_forward_nonvanishing}
For all $m\geq 1$, there exist constants $C>0$, $\lambda_0\geq 1$ and $\mu_0\geq 1$ such that for any $q_0\in L^2(\Omega,\mathcal F_0;L^2(\dom))$ and $G_i\in L^2_{\mathcal F}(0,T;L^2(\dom))$, $i=1,2$, the solution $y\in \mathcal W_T$ to \eqref{eq:forw_gen} satisfies 
\begin{align*} 
\esp&\left(\int_{Q_T}\theta^2\lambda\mu^2\xi|\nabla q|^2\dx\dt\right)+\esp\left( \int_{Q_T}\theta^2 \lambda^{3}\mu^{4}\xi^{3}|q|^2\dx\dt\right) + \esp\left(\int_{\Omega}\theta^{2}(T) \lambda^2 {\mu^{2}} |q(T)|^2   \dx \right)  \\
&\quad  \leq C\esp\left(\int_{Q_T}\theta^2|G_1|^2\dx\dt+\int_{Q_T}\theta^2\lambda^2\mu^2\xi^2|G_2|^2\dx\dt+\iint_{\mathcal \dom_0\times(0,T)}\theta^2\lambda^3\mu^4\xi^{3}|q|^2\dx\dt \right),
\end{align*}
for all $\mu\geq \mu_0$ and $\lambda\geq \lambda_0$.
\end{prop}

The proof of \Cref{prop:car_forward_nonvanishing} can be achieved by following the proof of \cite[Theorem 1.1]{liu14} with a few straightforward adaptations. For completeness, we give a brief sketch below. 

The starting point is to use \Cref{teo:contr_backward_source} with $F=\theta^2 \lambda^3\mu^4\xi^3 q$ and $y_T=-s^2\mu^{2}\theta^2(T)q(T)$, where $q$ is the solution to \eqref{eq:forw_gen} with given  $G_1$ and $G_2$. Observe that the weight functions in these data are well defined and bounded. We also remark that since the solution $q$ belongs to $ \mathcal W_T$, we have that $y_T=-\lambda^2\mu^2\theta^2(T)q(T)\in L^2(\Omega,\mathcal F_T;L^2(\dom))$ and thus system \eqref{eq:backward_source} with these given data is well-posed. 

Thus, from \Cref{teo:contr_backward_source}, we get that there exists a control $\widehat{h}\in L^2_{\mathcal F}(0,T;L^2(\dom))$ such that the solution $\widehat{y}$ to 
\begin{equation}\label{eq:backward_source_q}
\begin{cases}
\d\widehat{y}=(-\Delta \widehat{y}+\chi_{\dom_0}\widehat{h}+\theta^2\lambda^3\mu^4\xi^3 q)\dt+\widehat{Y}\d{W}(t) &\text{in }Q_T, \\
\widehat{y}=0 &\text{on }\Sigma_T, \\
\widehat{y}(T)=-\theta^2\lambda^2\mu^{2}q(T) &\text{in }\dom. 
\end{cases}
\end{equation}
satisfies $\widehat{y}(0)=0$ in $\dom$, a.s. Moreover, the following estimate holds
\begin{align}\notag 
\esp&\left(\int_{Q_T}\theta^{-2}|\widehat{y}|^2\dx\dt\right)+\esp\left(\int_{Q_T}\theta^{-2}\lambda^{-2}\mu^{-2}\xi^{-2}|\widehat{Y}|^2\dx\dt\right)+\esp\left(\int_{Q_T}\theta^{-2}\lambda^{-3}\mu^{-4}\xi^{-3}|\widehat{h}|^2\dx\dt\right) \\ \label{eq:est_for_carleman_weighted}
&\quad \leq C\esp\left(\int_{\dom}\lambda^2\mu^2\theta^2(T)|q(T)|^2\dx+\int_{Q_T}e^{-2s\varphi}s^{3}\xi^{3}|q|^2\dx\dt\right),
\end{align}
for some constant $C>0$ only depending on $\dom,\dom_0$. 

From \eqref{eq:forw_gen}, \eqref{eq:backward_source_q} and It\^{o}'s formula, we get
\begin{align*}
\esp&\left(\int_{\dom}\theta^2(T)\lambda^2\mu^{2}|q(T)|^2\dx\right)+\esp\left(\int_{Q_T}e^{-2s\varphi}s^3\xi^3|q|^2\dx\dt\right)\\
&=-\esp\left(\int_{Q_T}\widehat{y}G_1\dx\dt\right)-\esp\left(\int_{Q_T}\widehat{Y}G_2\dx\dt\right)-\esp\left(\int_0^{T}\!\!\!\int_{\dom_0}\widehat{h}q\dx\dt\right).
\end{align*}
Using Cauchy-Schwarz and Young inequalities, together with \eqref{eq:est_for_carleman_weighted}, it can be obtained from the above identity that
\begin{align*}
\esp&\left(\int_{\dom}\theta^2(T)\lambda^2\mu^2|q(T)|^2\dx\right)+\esp\left(\int_{Q_T}\theta^2\lambda^3\mu^4\xi^3|q|^2\dx\dt\right)\\
&\leq C\esp\left(\int_{Q_T}\theta^2|G_1|^2\dx\dt+\int_{Q_T}\theta^2\lambda^2\mu^2\xi^2|G_2|^2\dx\dt+\int_{0}^{T}\!\!\!\int_{\dom_0}\theta^2\lambda^3\mu^4\xi^3|q|^2\dx\dt\right).
\end{align*}
To add the integral containing $\nabla q$, it is enough to compute $\d(e^{-2s\varphi}\lambda\xi q^2)$ and argue as in Step 2 of the proof of \Cref{teo:contr_backward_source}. For brevity, we omit the details. 

\subsection{Globally Lipschitz nonlinearities depending on the gradient state}

It should be interesting to extend \Cref{th:semilinear_forward} to the case where the semilinearities $f$ and $g$ depend on the gradient of the state. More precisely, let us consider \begin{equation}\label{eq:forward_semilinear_grad}
\begin{cases}
\d{y}=(\Delta y + f(\omega,t,x,y, \nabla y)+\chi_{\dom_0}h)\dt+(g(\omega,t,x,y,\nabla y)+H)\d W(t) &\text{in }Q_T, \\
y=0 &\text{on }\Sigma_T, \\
y(0)=y_0 &\text{in }\dom,
\end{cases}
\end{equation}
where $f$ and $g$ are two globally Lipschitz nonlinear functions. We may wonder if \eqref{eq:forward_semilinear_grad} is small-time globally null-controllable. A good starting point seems to obtain a Carleman estimate for the backward equation \eqref{eq:system_z}, with a source term $\Xi \in L_{\mathcal{F}}^2(0,T;H^{-1}(\dom))$. This seems to be possible, according to \cite[Remark 1.4]{liu14}. By a duality argument, this would lead to a null-controllability result for the system \eqref{eq:sys_forward_source} similar to \Cref{teo:contr_forward_source}, with an estimate of $ \rho \widehat{y}$ in $L_{\mathcal F}^2(0,T;H_0^1(\dom))$, where $\rho$ is some suitable weight function. Details remain to be written.

\subsection{Extension of the method to other equations}

The method introduced in this article could probably be applied to other nonlinear equations for which there is a lack of compactness embeddings for the solutions spaces and for which we are able to derive Carleman estimates in the spirit of \cite[Theorem 2.5]{BEG16}. For instance, for the Schrödinger equation, it is a well known fact that there is no regularizing effect so there is a lack of compactness. Up to our knowledge, the following question is still open. Let $f : \C \rightarrow \C$ be a globally Lipschitz nonlinearity and $(T, \dom, \dom_0)$ be such that the so-called Geometric Control Condition holds. Is the system
\begin{equation*}\label{eq:schrodinger}
\begin{cases}
i \partial_t y =\Delta y + f(y) +\chi_{\dom_0}h&\text{in }Q_T, \\
y=0 &\text{on }\Sigma_T, \\
y(0)=y_0 &\text{in }\dom,
\end{cases}
\end{equation*}
globally null-controllable? See \cite{Zua03} or \cite{Lau14} for an introduction to this problem. We also refer to \cite{Lu13} for results in the stochastic setting.

\renewcommand{\abstractname}{Acknowledgements}
\begin{abstract}
\end{abstract}
\vspace{-0.5cm}
The work of the first author was supported by the programme ``Estancias posdoctorales por M\'exico'' of CONACyT, Mexico. The work of the second author has been supported by the SysNum cluster of excellence University of Bordeaux.

\bibliographystyle{alpha}
\small{\bibliography{bib_nonlin_stoch}}

\end{document}